\documentclass[a4paper,12pt]{article}

\usepackage{float}
\usepackage{amsmath} 
\usepackage{amssymb}
\usepackage{amsfonts}

\usepackage{epstopdf}
\usepackage{color,fancybox,graphicx}
\usepackage{url}
\usepackage{psfrag}
\usepackage{color}
\usepackage{algorithm, algorithmic}
\usepackage{pstricks}
\usepackage{graphicx,psfrag}
\usepackage{epsfig}

\numberwithin{equation}{section}

\newtheorem{lem}{Lemma}[section]
\newtheorem{cor}{Corollary}[section]
\newtheorem{pro}{Proposition}[section]
\newtheorem{theo}{Theorem}[section]

\newcommand{\1}{\mathbf{1}}

\newcommand{\un}{{\mathbf{1}}}
\newcommand{\A}{\mathbf{A}}
\newcommand{\B}{\mathbf{B}}
\newcommand{\F}{{\cal F}}
\renewcommand{\d}{\delta}

\newcommand{\g}{\gamma}
\newcommand{\e}{\eta}
\newcommand{\s}{\sigma}
\renewcommand{\a}{\alpha}
\newcommand{\V}{\mathbb V}

\renewcommand{\div}{\mathrm{div}}
\def\R{\mathbb{R}}
\def\P{\mathbb{P}}
\def\E{\mathbb{E}}

\newcommand{\ck}{ {\bf k }}
\newcommand{\ccr}{ {\bf r }}

\newcommand{\Ba}{ {\cal B }}

\newcommand{\point}{\mbox{\LARGE .}}

\newcommand{\dbar}{\mathchar'26\mkern-12mu d}

\begin{document}

\begin{center}
{\sc \Large Fluctuation Analysis
 \\
\vspace{0.2cm}
of Adaptive Multilevel Splitting}
\vspace{0.5cm}

\end{center}

{\bf Fr\'ed\'eric C\'erou\footnote{Corresponding author.}}\\
{\it INRIA Rennes \& IRMAR, France }\\
\textsf{frederic.cerou@inria.fr}
\bigskip


{\bf Arnaud Guyader}\\
{\it Universit\'e Pierre et Marie Curie, France }\\
\textsf{arnaud.guyader@upmc.fr}

\medskip

\begin{abstract}
\noindent {\rm Multilevel Splitting, also called Subset Simulation, is a Sequential Monte Carlo method to simulate realisations of a rare event as well as to estimate its probability. This article is concerned with the convergence and the fluctuation analysis of Adaptive Multilevel Splitting techniques. In contrast to their fixed level version, adaptive techniques estimate the sequence of levels on the fly and in an optimal way, with only a low additional computational cost. However, very few convergence results are available for this class of adaptive branching models, mainly because the sequence of levels depends on the occupation measures of the particle systems. This article proves the consistency of these methods as well as a central limit theorem. In particular, we show that the precision of the adaptive version is the same as the one of the fixed-levels version where the levels would have been placed in an optimal manner.
\medskip

\medskip

}

\end{abstract}

\section{Introduction}
Multilevel Splitting techniques were introduced as natural heuristics in the 1950s' by Kahn and Harris \cite{kahn} and Rosenbluth and Rosenbluth \cite{dim21} to analyze particle transmission energies and molecular polymer conformations. In their basic form, these methods can be interpreted as a genetic mutation-selection algorithm. The mutation transition reflects the free evolution of the physical model at hand, while the selection is an acceptance-rejection transition equipped with a recycling mechanism. The samples entering a critical level are more likely to be selected and duplicated. The  genealogy associated with these genetic type particles model represents the statistical behavior of the system passing through a cascade of critical rare events. 
\medskip

Interestingly, these models can also be seen as a mean field particle approximation of Feynman-Kac measures. This interpretation depends on the application area of interest. In scientific computing and mathematical biology, these stochastic techniques are often termed genetic algorithms. In machine learning and advanced signal processing, they are referred as Sequential Monte Carlo or Particle Filters. In computational and quantum physics, they belong to the class of Diffusion Monte Carlo methods. The analysis of this class of branching and mean field type particle methods is now well understood (see, for instance, ~\cite{cdg,alea06,del2013mean,dm-2000} and references therein).\medskip

The present article is concerned with the convergence analysis of a more sophisticated class of adaptive particle methods where both the selection functions and the mutation transitions depend on the occupation of the system. The selection functions are chosen to sequentially achieve a prescribed proportion of samples in an higher critical level set, while the mutation transitions are dictated by some Markov Chain Monte Carlo (MCMC) methods on the adaptive critical level sets. This adaptive multilevel technique is a natural and popular approach amongst practitioners, but there are very few convergence results for this class of models. 
\medskip

The first analysis of this class of models has been developed by Del Moral, Doucet and Jasra in~\cite{ddj-02}, in terms of adaptive resampling times associated with some criteria such as the effective sample size. Nonetheless, their result only applies to adaptive models associated with parametric level sets and equipped with sequential resampling times related to some fixed threshold. Thus, this does not correspond precisely to the purpose of Adaptive Multilevel Splitting methods that we are presently interested in. 
\medskip

In a slightly different framework, the recent article by Beskos, Jasra, Kantas and Thi\'ery \cite{bjkt14} is also related to the present paper. Specifically, the authors present a detailed analysis for a class of adaptive Sequential Monte Carlo models under regularity properties on the dependency of the mutation transitions and the selection functions with respect to the occupation measures of the system. The proofs in \cite{bjkt14} reveal that these regularity properties are essential to develop a first order perturbation analysis between the adaptive particle models and their limiting measures. Unfortunately, this framework does not apply to indicator selection functions arising in classical multilevel splitting methodologies and developed in the present article. As a consequence, even if the global goal here is roughly the same as in \cite{ddj-02,bjkt14}, the techniques developed for establishing our convergence results are quite different. Note also that in the context of adaptive tempering (a context considered  in \cite{bjkt14}), Giraud and Del Moral give non-asymptotic bounds on the error in \cite{gdelmo}. 
\medskip

Let us first specify our framework and notation. In all the paper, we suppose that $X$ is a random vector in $\R^d$ with law $\eta$ that we can simulate, and $S$ is a mapping from $\R^d$ to $\R$, also called a score function. Then, given a threshold $L^\star$ which lies far out in the right hand tail of the distribution of $S(X)$, our goal is to estimate the rare event probability $P={\mathbb P}(S(X)>L^\star)$. This very general context includes applications in queuing networks, insurance risks, random graphs (as found in social networks, or epidemiology), etc., see e.g. \cite{gobetliu} for some of them, and a discussion on practical implementations.\medskip

In this context, a crude Monte Carlo uses an i.i.d. $N$-sample $X_1,\ldots,X_N$ to estimate $P$ by the fraction $\hat P_{mc} = \#\{i:\ S(X_i) > L^\star\}/N$. However, in order to obtain a reasonable precision of the estimate given by the relative variance $\V(\hat P_{mc})/P^2=(1-P)/(NP)$, one needs a sample size $N$ of order at least $P^{-1}$. Obviously, this becomes unrealistic when $P$ is very small, hence the use of variance reduction techniques. 
\medskip

Importance Sampling, which draws samples according to $\pi$ and weights each
observation $X=x$ by $w(x)={d\eta(x)}/{d\pi(x)}$, may decrease the variance of the estimated probability dramatically, which in turn reduces the need for such large sample sizes. We refer to Robert and Casella \cite{robert+2004} for a discussion on Importance Sampling techniques in general, and to Bucklew \cite{bucklew04} and L'Ecuyer, Mandjes and Tuffin \cite[Chapter 2]{rubino} for the application in the context of rare event estimation. Notice that, in rare event estimation, it is customary to design an importance sampling scheme using a large deviation principle. Although it often gives an efficient method, this approach may fail dramatically, even compared to crude Monte Carlo, when the rare event has two or more most likely occurrences. As explained by Glasserman and Wang  in the introduction of \cite{GW}, ``Simply put, an analysis of a first moment cannot be expected to carry a guarantee about the behavior of a second moment.''
\medskip

Multilevel Splitting represents another powerful algorithm for rare event estimation. The basic idea of Multilevel Splitting, adapted to our problem, is to fix a set of increasing levels $-\infty = L_{-1}<L_0 < \dots < L_{n-1}<L_n=L^\star$, and to decompose the tail probability thanks to Bayes formula, that is
$$\P(S(X) > L^\star) = \prod_{p=0}^{n} \P(S(X) > L_{p} | S(X) > L_{p-1}).$$
Each conditional probability $\P(S(X) > L_{p} | S(X) > L_{p-1})$ is then estimated separately. We refer the reader to L'Ecuyer, Le Gland, Lezaud and Tuffin \cite[Chapter 3]{rubino} for an in-depth review of the Multilevel Splitting method and a detailed list of references. Two practical issues associated with the implementation of Multilevel Splitting are: first, the need for computationally efficient algorithms for estimating the successive conditional probabilities; second, the optimal selection of the sequence of levels.
\medskip  

The first question can be addressed thanks to the introduction of Markov Chain Monte Carlo procedures at each step of the algorithm. This trick was proposed in different contexts and through slightly different variants by Au and Beck \cite{Au2001263,au:901}, Del Moral, Doucet and Jasra \cite{delmoral06b}, Botev and Kroese~\cite{kroese08a}, Rubinstein~\cite{rubinstein08a}. 
\medskip

The second question is straightforward in the idealized situation where one could estimate the successive quantities $\P(S(X) > L_{p}|S(X) > L_{p-1})$ independently at each step. Indeed, considering the variance of the estimator, it is readily seen that the best thing to do is to place the levels as evenly as possible in terms of the intermediate probabilities, that is to take, for all $p$, 
$$\P(S(X) > L_{p}|S(X) > L_{p-1})=\P(S(X) > L^\star)^{\frac{1}{n+1}}.$$ 
But, since little might be known about the mapping $S$, the only way to achieve this goal is to do it on the fly by taking advantage of the information of the current sample at each step. This method is called Subset Simulation (see Au and Beck \cite{Au2001263,au:901}) or Adaptive Multilevel Splitting (see C\'erou and Guyader \cite{cg2}), and may be seen as an adaptive Sequential Monte Carlo method specifically dedicated to rare event estimation.
\medskip

However, except in the idealized situation where one considers a new independent sample at each step (see C\'erou, Del Moral, Furon and Guyader \cite{cdfg}, Guyader, Hengartner and Matzner-L\o ber \cite{ghm}, Br\'ehier, Leli\`evre and Rousset \cite{blr14}, and Simonnet \cite{s14}), there are only very few results about the theoretical properties of this efficient algorithm. From a broader point of view, as duly noticed in ~\cite{ddj-02,bjkt14}, this disparity between theory and practice holds true for adaptive Sequential Monte Carlo methods in general. As such, the present article is in the same vein as \cite{ddj-02,bjkt14} and might be seen as a new step towards a better understanding of the statistical properties of adaptive Sequential Monte Carlo methods. 
\medskip

In particular, the take-home message here is the same as in \cite{ddj-02,bjkt14}, namely that the asymptotic variance of the adaptive version is the same as the one of the fixed-levels version where the levels would have been placed in an optimal manner. However, there are substantial differences between  \cite{ddj-02,bjkt14} and the present contribution.\medskip

In \cite{ddj-02}, the adaptive parameter is the time at which one needs to resample. This approach can be used for rare event if we choose a possibly long sequence of deterministic  levels $L_1,\dots,L_n$, and resample only when the current level sees a given proportion of particles to be already killed. The authors provide convergence results, including a CLT, when the number $N$ of particles goes to $\infty$, but for fixed levels $L_1,\dots,L_n$. To get the kind of results of the present contribution, one would need to let also $n$ go to $\infty$, and this cannot be achieved by the coupling technique used in \cite{ddj-02} due to the inherent  jittering of the adaptive levels, which is typically of order $1/\sqrt{N}$. If the granularity of the levels goes to $0$ as $N$ goes to $\infty$, then there is little hope  that the adaptive particle system coincides with the optimal one with large probability as in their Theorem 2.3.\medskip 

In  \cite{bjkt14}, the authors consider different scenarios, including adaptive proposal and adaptive tempering, where they can make a Taylor expansion of the adaptive selection function, and the adaptive kernel, in the vicinity of the optimal parameter. This leads to additional terms in the asymptotic variance that may cancel in some cases (adaptive proposal), giving the same variance as in the non adaptive optimal case. Yet, let us emphasize again that the inherent unsmoothness of the selection functions of interest here (going abruptly from $0$ to $1$ when crossing a level set for $S$)  leads to different proofs, meaning that their results and even techniques, although very interesting in and by themselves, can definitely not be applied in our context.\medskip

The paper is organized as follows. In Section \ref{framework}, we introduce some notation and describe the Multilevel Splitting algorithms. The asymptotic results (laws of large numbers and central limit theorems) are presented in Section \ref{results}. Section \ref{assumption} comes back on the assumption required for our CLT type result to be valid. Section \ref{azp} is devoted to the proofs of the theorems, while technical results are postponed to Section \ref{technical}.

\section{Multilevel splitting techniques}\label{framework}
\subsection{Framework and notation }
We consider an $\R^d$-valued random variable $X$ with distribution $\eta$, for some $d\geq 1$. We assume that $\eta$ has a density with respect to Lebesgue's measure $dx$ on $\R^d$ and, by a slight abuse of notation, we denote $\eta(x)$ this density. We also consider a mapping $S$ from $\R^d$ to $\R$. If $S$ is Lipschitz with $|DS|>0$ almost everywhere, where $|DS|$ stands for the Euclidean norm of the gradient of $S$, then the coarea formula (see for example \cite{evans}, page 118, Proposition 3) ensures that the random variable $Y=S(X)$ is absolutely continuous with respect to Lebesgue's measure on $\R$, and its density is given by the formula
\begin{equation}\label{fy}
f_Y(s)=\int_{S(x)=s}\eta(x)\frac{\dbar x}{|DS(x)|},
\end{equation}
where $\dbar x$ stands for the Hausdorff measure on the level set $S^{-1}(s)=\{x\in\R^d, S(x)=s\}$. In this notation, given $\alpha\in(0,1)$, the $(1-\alpha)$ quantile of $Y$ is simply $F_Y^{-1}(1-\alpha)$, where $F_Y$ stands for the cumulative distribution function (cdf for short) of $Y$.
\medskip

Consider a real number (or level) $L^\star$ lying far away in the right hand tail of $S(X)$ so that the probability $P=\P(Y\geq L^\star)$ is very small. For any bounded and measurable function $f:\R^d\to\R$ (denoted $f\in\Ba(\R^d)$ in all the paper) which is null below $L^\star$ (implicitly: with respect to $S$), our goal is to estimate its expectation with respect to $\e$, that is the quantity
\begin{equation}\label{esp}  
E=\E[f(X)]=\E[f(X)\un_{S(X)\geq L^\star}].
\end{equation}

To this end, we fix an $\a\in(0,1)$ (in practice one may typically choose $\a=3/4$), and consider the decomposition
\begin{equation}\label{p}  
P=\P(Y\geq L^\star)=r\times \a^{n}\quad\mbox{with}\quad n=\left\lfloor\frac{\log \P(Y\geq L^\star)}{\log\a}\right\rfloor,
\end{equation}   
so that $r\in(\a,1]$. For the sake of simplicity and since this is always the case in practice, we assume that $r$ belongs to the open interval $(\alpha,1)$. With the convention $L_{-1}=-\infty$, we define the increasing sequence of levels $(L_p)_{p\geq -1}$ as follows
$$L_0=F_Y^{-1}(1-\a)<\dots<L_{n-1}=F_Y^{-1}(1-\a^n)<L^\star<L_{n}=F_Y^{-1}(1-\a^{n+1}).$$

Once and for all, we assume that the density $f_Y$, as defined in equation (\ref{fy}), is {\bf continuous} and {\bf strictly positive} at each $L_p$, for $p\in\{0,\dots,n\}$. This will guarantee that the quantiles are well defined and that the empirical ones have good convergence properties. 
\medskip

Following the notations of \cite{delmoral04a,del2013mean}, we associate to these successive levels the potential functions 
$$\forall -1\leq p<n,\qquad G_p=\un_{{\cal A}_{p}}\quad\mbox{\rm with}\quad {\cal A}_{p}=\{x\in\R^d: S(x)\geq L_{p}\}.$$
The restriction of $\e$ to ${\cal A}_{p-1}$ is then denoted $\e_{p}$. More formally, we have
$$\e_p(dx)=\a^{-p}\un_{{\cal A}_{p-1}}(x)\eta(x)dx=\a^{-p}G_{p-1}(x)\eta(x)dx.$$ 
By construction, we have 
$$\eta_p(G_p)=\eta_p(\un_{{\cal A}_p})=\P\left(S(X)\geq L_{p}|S(X)\geq L_{p-1}\right)=\alpha.$$
We also notice that the interpolating measures $\e_p$ are connected by the Boltzmann-Gibbs transformation
$$\eta_{p+1}(dx)=\Psi_{G_p}(\eta_p)(dx)=\frac{1}{\eta_p(G_p)}G_p(x)\eta_p(dx)=\alpha^{-1}G_p(x)\eta_p(dx).$$
Moreover, we consider a collection of Markov transitions from ${\cal A}_{p-1}$ into itself defined for any $x\in{\cal A}_{p-1}$ by
\begin{equation}\nonumber
M_p(x,dx')=K_p(x,dx')\un_{{\cal A}_{p-1}}(x')+K_p(x,\bar{{\cal A}}_{p-1})~\d_x(dx'),
\end{equation}
where $\bar{{\cal A}}_{p-1}=\R^d-{\cal A}_{p-1}$, and $K_p$ stands for a collection of $\eta$-reversible Markov transitions on $\R^d$, meaning that for all $p$ and all couple $(x,x')$, we have the detailed balance equation
\begin{equation}\label{dbe}
\eta(dx)K_p(x,dx')=\eta(dx')K_p(x',dx).
\end{equation}
We extend $M_p$ into a transition kernel on $\R^d$ by setting $M_p(x,dx')=\delta_x(dx')$ whenever $x\not\in {\cal A}_{p-1}$. Under the assumption that $K_p$ is $\e$-symmetric, it is easy to check that $M_p$ is $\e_p$-invariant, meaning that $\e_pM_p=\e_p$ for all $p\geq1$. In addition, we have the recursion
$$\e_p(dx^{\prime})=\a^{-1}(\e_{p-1}Q_p)(dx^{\prime})=\a^{-1}\int \e_{p-1}(dx)Q_p(x,dx^{\prime}),$$
with the integral operators
$$Q_p(x,dx^{\prime})=G_{p-1}(x)M_p(x,dx^{\prime}).$$
Next, let us denote $(X_p)_{p\geq 0}$ a non homogeneous Markov chain with initial distribution $\e_0=\e$ and elementary transitions $M_{p+1}$. In this situation, it is readily seen that
\begin{equation}\label{ljnx}
\alpha^n~\e_n(f)=\E\left[f(X_n)\prod_{q=0}^{n-1}G_{q}(X_q)\right]\quad\Longleftrightarrow\quad \alpha^n~\e_n=\e_0Q_{0,n}.
\end{equation}
with the Feynman-Kac semigroup $Q_{0,n}$ associated with the integral operators $Q_p$ defined by $$\forall 0\leq p\leq n\qquad Q_{p,n}=Q_{p+1}Q_{p+1,n}$$
In this notation, we have
\begin{eqnarray*}
E&=&\E[f(X)]=\E[f(X)\un_{S(X)\geq L^\star}]=\a^n\times\e_n(f\times\un_{S(\cdot)\geq L^\star})\\
P&=&\P(Y\geq L^\star)=\P(S(X)\geq L^\star)=\a^{n}\times\eta_n(\un_{S(\cdot)\geq L})=\a^{n}\times r
\end{eqnarray*}
and
\begin{equation}\label{C}
f=f\times\un_{S(\cdot)\geq L^\star}\ \Longrightarrow\ C=\E[f(X)|S(X)\geq L^\star]=\frac{\e_{n}(f)}{\e_n(\un_{S(\cdot)\geq L^\star})}=\frac{\e_{n}(f)}{r}.
\end{equation}

We will now describe two multilevel splitting techniques in order to estimate these quantities. The optimal Feynman-Kac particle approximation of the flow (\ref{ljnx}) corresponds to the fixed-levels method that we describe in Section \ref{fx}. As this approximation is not possible in practice, we detail in Section \ref{ad} the corresponding adaptive Feynman-Kac particle approximation, known as Adaptive Multilevel Splitting or Subset Simulation.  

\subsection{The fixed-levels method}\label{fx}
Following the notation of \cite{delmoral04a}, the fixed-levels approximation of the flow (\ref{ljnx}) works as follows. Let $(X_p^1,\dots,X_p^N)_{0\leq p\leq n}$ be an $(\R^d)^N$-valued Markov chain with initial distribution $\e_0^{\otimes N}$ and for which each elementary transition $X_{p}^i\leadsto X_{p+1}^i$ is decomposed into the following separate mechanisms:
\begin{enumerate}
\item Selection step: compute $\check{\e}^N_{p}(G_p)$, which is the proportion of the sample $(X_p^1,\dots,X_p^N)$ such that $S(X_p^i)\geq L_{p}$.
\item Multinomial step: from the $\check{\e}^N_{p}(G_p)N$-sample with distribution $\eta_{p+1}$, draw an $N$-sample $(X_{p+1/2}^1,\dots,X_{p+1/2}^N)$ with the same distribution.
\item Transition step: each $X_{p+1/2}^i$ evolves independently to a new site $X_{p+1}^i$ randomly chosen with distribution $M_{p+1}(X_{p+1/2}^i,dx')$.
\item Incrementation step: $p=p+1$. If $p=n$, then stop the algorithm, else go to step 1 (selection step). 
\end{enumerate}

Let us denote $\check{\g}_n^N(1)$ the normalizing constant defined by
$$\check{\g}_n^N(1)=\prod_{p=0}^{n-1}\check{\e}^N_{p}(G_p).$$
In our framework, its deterministic counterpart is simply 
$$\g_n(1)=\prod_{p=0}^{n-1}\e_{p}(G_p)=\a^n.$$
For any $f\in\Ba(\R^d)$, the normalized and unnormalized measures $\check{\e}_{n}^N(f)$ and $\check{\g}_n^N(f)$ are respectively defined by
$$\check{\e}_{n}^N(f)=\frac{1}{N}\sum_{i=1}^Nf(X_n^i)\quad\mbox{and}\quad\check{\g}_n^N(f)=\check{\g}_n^N(1)\times\check{\e}_{n}^N(f).$$
The fixed-levels algorithm provides the following estimates:
\begin{enumerate}
\item[$(i)$] The estimate of the expectation $E=\E[f(X)\un_{S(X)\geq L^\star}]=\g_{n}(f\times\un_{S(\cdot)\geq L^\star})$ is given by $\check{E}=\check{\g}_n^N(f\times\un_{S(\cdot)\geq L^\star})$. 
\item[$(ii)$] The rare event probability $P=\P(S(X)\geq L^\star)$ is estimated by the quantity $\check{P}=\check{\g}_n^N(\un_{S(\cdot)\geq L^\star})$.
\item[$(iii)$] The estimate of the conditional expectation $C=\E[f(X)|S(X)\geq L^\star]$ is 
$$\check{C}=\frac{\check{\e}_n^N(f\times\un_{S(\cdot)\geq L^\star})}{\check{\e}_n^N(\un_{S(\cdot)\geq L^\star})}=\frac{\sum_{i=1}^Nf(X_{n}^i)\un_{S(X_n^i)\geq L^\star}}{\sum_{i=1}^N\un_{S(X_n^i)\geq L^\star}}.$$ 
\end{enumerate}

These particle models associated with a collection of deterministic potential functions $G_p$ and Markov transitions $M_p$ belong to the class of Feynman-Kac particle models. This class of mean field particle models has been extensively studied in a very general context, including the asymptotic behavior as the number $N$ of particles goes to infinity. We refer the reader to~\cite{delmoral04a} and the more recent research monograph~\cite{del2013mean}, with references therein. We will recall some of these results in Section \ref{comparison}.\medskip

In our specific context, the obvious drawback of these Feynman-Kac particle approximations is the impossibility to fix in advance the successive levels $L_0,\dots,L_n$, hence the use of adaptive methods that we describe in the following section.

\subsection{The adaptive method}\label{ad} 
An efficient way to estimate the quantities $E$, $P$ and $C$ is to use Adaptive Multilevel Splitting methods. To describe with some precision these particle splitting models, it is convenient to consider a collection of potential functions and Markov transitions indexed by $\R$. Thus, for any real number $L$, we set 
$$G_L=\un_{{\cal A}_L}\quad \mbox{with}\quad {\cal A}_L=\{x\in\R^d: S(x)\geq L\}.$$
We also consider the collection of Markov transitions from ${\cal A}_{L}$ into itself defined for any $x\in{\cal A}_{L}$ by
$$M_{p,L}(x,dx')=K_p(x,dx')\un_{{\cal A}_L}(x')+K_p(x,\bar{{\cal A}}_L)\d_x(dx').$$
As before, we extend $M_{p,L}$ into a transition kernel on $\R^d$ by setting $M_{p,L}(x,dx')=\delta_x(dx')$ whenever $x\not\in {\cal A}_L$, and we set
$$Q_{p,L}(x,dx^{\prime})=G_{L}(x)M_{p,L}(x,dx^{\prime}).$$
In this slight abuse of notation, we have 
$$
L=L_{p-1}\quad\Longrightarrow\quad (G_L,{\cal A}_L)=(G_{p-1},{\cal A}_{p-1})\quad\mbox{\rm and}\quad (M_{p,L},Q_{p,L})=(M_p,Q_p).
$$
Of special interest will be the case where $L$ is a given quantile. We distinguish two cases: 
\begin{itemize}
\item Firstly, for any positive and 
 finite measure $\nu$ on $\R^d$ with a density  with respect to Lebesgue's measure, the level $L_\nu$ is defined as the $(1-\alpha)$ quantile of the probability measure $(S_\ast \nu)/\nu(\R^d)$, that is
\begin{equation}\label{lknc}
L_\nu=L(\nu)=F_\nu^{-1}(1-\a)\quad \mbox{where}\quad F_\nu(y)=\frac{\nu(S^{-1}((-\infty,y]))}{\nu(\R^d)}.
\end{equation}
In order to lighten the notations a bit, we will write 
$$G_\nu:=G_{L_\nu}\qquad {\cal A}_\nu:={\cal A}_{L_\nu}\ \qquad M_{p,\nu}:=M_{p,L_\nu}$$ 
and\label{ppu1}
$$Q_{p,\nu}(x,dx^{\prime})=G_\nu(x)M_{p,\nu}(x,dx^{\prime}).$$
\item Secondly, given a sample of vectors $(X_i)_{1\leq i\leq N}$ in $\R^d$, we consider an auxiliary sequence of i.i.d. uniformly distributed random variables $(U_1,\dots,U_N)$ and the following total order on the couples $(X_i,U_i)_{1\leq i\leq N}$:
\begin{align}
&(X_i,U_i)<(X_j,U_j)\ \Leftrightarrow\ S(X_i)<S(X_j)\ \mbox{or}\ S(X_i)=S(X_j)\ \mbox{and}\ U_i<U_j.\label{ordre}
\end{align}
Obviously, since the $U_i$'s are uniformly distributed, equality between two couples almost surely never happens. Hence we can consider the associated order statistics
\begin{align}
&(X_{(1)},U_{(1)})<\dots<(X_{(n)},U_{(n)}),\nonumber
\end{align}
and we define the empirical $(1-\alpha)$ quantile $L^N$ as 
\begin{align}
&L^N:=(S(X_{(\lfloor N(1-\alpha)\rfloor)}),U_{(\lfloor N(1-\alpha)\rfloor)}).\label{defquant}
\end{align} 
In particular, one can notice that the number of couples strictly above $L^N$ (with respect to the previous order) is equal to $\lceil N\alpha\rceil$. 
\end{itemize}

{\bf Remarks:} 
\begin{itemize}
\item In order to lighten the writings, we will usually not mention that some auxiliary uniform random variables $(U_1,\dots,U_N)$ are {\bf always} attached to a sample $(X_1,\dots,X_N)$. In particular, this will be implicit to define the associated empirical quantile $L^N$, and the relation $S(X_i)>L^N$ must be understood accordingly. Otherwise, the notation $L^N$ refers only to its first component in definition (\ref{defquant}). However, considering the context, there should be no ambiguity. Note also that when considering the convergence of empirical quantiles, only the first component will be considered, as there is no reason why the uniform random variable would converge.
\item In our situation, it turns out that equality between several $X_i$'s\ will have no influence on the CLT type result we want to establish. Indeed, ties come from the multinomial step of the algorithm below, but as mentioned in the proof of Lemma \ref{lemmulti}, one can control these events very precisely (see for example Chapter 4 of \cite{lucbucket}). As we will see, the maximum number of particles on the same level set for $S$, at stage $q$, is typically $o_p((\log N)^q)$, making the mass at a single location  $o_p((\log N)^q/N)$, while the error between the conditional measure $\eta_q$ and its particle approximation will be of order $1/\sqrt{N}$ as expected.\medskip
\end{itemize}

In this context, the {\bf adaptive} particle approximation of the flow (\ref{ljnx}) is defined in terms of
an $(\R^d)^N$-valued Markov chain $(X_p^1,\dots,X_p^N)_{p\geq 0}$ with initial distribution $\e_0^{\otimes N}$.  We start with $p=0$ and a $\e_0^{\otimes N}$ sample $(X_0^1,\dots,X_0^N)$. The  elementary transitions 
$X_{p}^i\leadsto X_{p+1}^i$ are decomposed into the following separate mechanisms:
\begin{enumerate}
\item Quantile step: compute the empirical $(1-\alpha)$ quantile $L^N_{p}$ of the sample $(X_p^1,\dots,X_p^N)$ in the sense of (\ref{defquant}). If $L^N_{p}\geq L^\star$, then stop the algorithm, else go to step \ref{2} (multinomial step).
\item Multinomial step: draw an $N$-sample $(X_{p+1/2}^1,\dots,X_{p+1/2}^N)$ with common distribution
$$\tilde{\eta}_p^N(dx)=\frac{1}{\lceil N\alpha\rceil}\sum_{i:~S(X_{p}^i)> L^N_{p}}~\delta_{X_p^i}(dx)=\frac{1}{\lceil N\alpha\rceil}\sum_{i=1}^{\lceil N\alpha\rceil}~\delta_{\tilde X_p^i}(dx).$$\label{2}
\item Exploration step: each $X_{p+1/2}^i$ evolves independently to a new site $X_{p+1}^i$ randomly chosen with distribution $M_{p+1,L^N_{p}}(X_{p+1/2}^i,dx')$. 
\item Incrementation step: $p=p+1$.  Go to step 1.
\end{enumerate}
Denote $\hat{n}$ the last index $p$ such that $L_p^N<L^\star$. This algorithm provides the following estimates:
\begin{enumerate}
\item[$(i)$] The estimate of the expectation $E=\E[f(X)]=\E[f(X)\un_{S(X)\geq L^\star}]$ considered in (\ref{esp}) is 
$$\hat{E}=\a^{\hat{n}}\times\e_{\hat{n}}^N(f)=\a^{\hat{n}}\times\frac{1}{N}\sum_{i=1}^Nf(X_{\hat{n}}^i)\un_{S(X_{\hat{n}}^i)\geq L^\star}.$$ 
\item[$(ii)$] The rare event probability $P=\P(S(X)\geq L^\star)$ considered in (\ref{p}) is estimated by the quantity
$$\hat{P}=\a^{\hat{n}}\times\eta_{\hat{n}}^N(\un_{S(\cdot)\geq L^\star})=\a^{\hat{n}}\times\frac{1}{N}\sum_{i=1}^N\un_{S(X_{\hat{n}}^i)\geq L^\star}.$$
\item[$(iii)$] For the conditional expectation $C=\E[f(X)|S(X)\geq L^\star]$ considered in (\ref{C}), still with $f=f\times\un_{S(\cdot)\geq L^\star}$, the estimate is 
$$\hat{C}=\frac{\e_{\hat{n}}^N(f)}{\e_{\hat{n}}^N(\un_{S(\cdot)\geq L^\star})}=\frac{\sum_{i=1}^Nf(X_{\hat{n}}^i)\un_{S(X_{\hat{n}}^i)\geq L^\star}}{\sum_{i=1}^N\un_{S(X_{\hat{n}}^i)\geq L^\star}}.$$ 
\end{enumerate}
The purpose of Section \ref{adres} is to expose some asymptotic results on these estimators. 

\subsection{Metropolis-Hastings kernels}\label{mh}
Let us briefly recall Metropolis-Hastings algorithm \cite{metropolis53,hastings1970}, which is a possible way to obtain a collection $K_p$ of $\eta$-reversible Markov transitions. We emphasize that, from a practical viewpoint, the kernels $K_p$ are a key ingredient of the previous algorithms, for fixed levels as well as for adaptive ones. Hereafter we follow the presentation of \cite{tierney94}.\medskip

Let $\ck_p$ be a Markov transition kernel of the form
$$\ck_p(x,dx')=\ck_p(x,x')dx'.$$ 
Let $E^+=\{x\in\R^d,\ \eta(x)>0\}$ and, for the sake of simplicity, assume that $\ck_p(x,E^+)=1$ for any $x\notin E^+$. Next define the acceptance ratio
$$a_p(x,x')=\left\{\begin{array}{ll}
\min\left\{\displaystyle\frac{\eta(x')\ck_p(x',x)}{\eta(x)\ck_p(x,x')},1\right\}&\ \mbox{if}\ \eta(x)\ck_p(x,x')>0,\\
1&\ \mbox{if}\ \eta(x)\ck_p(x,x')=0.\end{array}\right.$$  
The success of Metropolis-Hastings algorithm comes from the fact that it only depends on $\eta$ through ratios of the form $\eta(x')/\eta(x)$, hence $\eta$ only needs to be known up to a normalizing constant. If we define the off-diagonal density of a Metropolis kernel as 
$$\ck_p^a(x,x')=a_p(x,x')\ck_p(x,x'),$$
and set
$$\ccr_p^a(x)=1-\int_{\R^d}\ck_p^a(x,x')dx',$$
then the Metropolis kernel $K_p$ can be written as 
$$K_p(x,dx')=\ck_p^a(x,x')dx'+\ccr_p^a(x)\delta_x(dx').$$
Since $\ck_p^a$ satisfies the detailed balance equation
\begin{equation}\label{kpa}
\eta(x)\ck_p^a(x,x')=\eta(x')\ck_p^a(x',x),
\end{equation}
it follows that $K_p$ is an $\eta$-reversible Markov transition kernel. Finally, let us mention that, for any function $\varphi\in\Ba(\R^d)$, we will denote $\ck_p^a(\varphi)$ the function defined by
$$\ck_p^a(\varphi)(x)=\int_{\R^d}\ck_p^a(x,x')\varphi(x')dx'=\int_{\R^d}a_p(x,x')\ck_p(x,x')\varphi(x')dx',$$
so that
\begin{equation}\label{kpaa} 
K_p(\varphi)(x)=\ck_p^a(\varphi)(x)+\ccr_p^a(x)\times\varphi(x).
\end{equation}
This expression will be useful in the proof of Proposition \ref{VieillePuteAlbanaise}.

\section{Consistency and fluctuation analysis}\label{results}

\subsection{Adaptive Multilevel Splitting}\label{adres}
We prove in Theorem \ref{kjsc} the almost sure convergence of $L_p^N$ to $L_p$. As a byproduct, we deduce that the probability that the algorithm does not stop after the right number of steps (i.e., that $\hat{n}\neq n$) goes to zero when $N$ goes to infinity. Then, in Theorem \ref{theorem1}, we focus our attention on the fluctuations of $\e_n^N(f)$ around $\e_n(f)$.  

\begin{theo}\label{kjsc}
For all $p\in\{0,\dots,n\}$,
$$
L_p^N\xrightarrow[N\to\infty]{a.s.}L_{p}.
$$
Besides, for all $f\in L^2(\eta)$,
$$
\eta_p^N(f)\xrightarrow[N\to\infty]{\P}\eta_p(f),
$$
and for all $f\in\Ba(\R^d)$,
$$
\eta_p^N(f)\xrightarrow[N\to\infty]{a.s.}\eta_p(f).
$$
\end{theo}

Note that a consequence of Theorem \ref{kjsc} is that the couple $(L_{n-1}^N,L_{n}^N)$ converges almost surely to $(L_{n-1},L_{n})$. As claimed before, this ensures that, almost surely for $N$ large enough, $L_{n-1}^N<L^\star<L_{n}^N$, which means that $\hat{n}=n$.
\medskip
  
The fluctuations of $\e_n^N$ around the limiting measure $\e_n$ are expressed in terms of the normalized Feynman-Kac semigroups
$\overline{Q}_{q,p}$ defined by
$$\forall 0\leq q\leq p\leq n,\qquad
\overline{Q}_{q,p}=\frac{Q_{q,p}}{\eta_q(Q_{q,p}(1))}=\alpha^{q-p}\times Q_{q,p}.$$

We also need to specify some regularity assumptions on the score function $S$ and the transition kernels $K_q$ for which our CLT type result is valid. For any $q>0$, we first introduce the set of functions 
$${\cal B}_q=\left\{g:\R^d\to\R,\ \exists (g_0\dots,g_{q-1})\in\Ba(\R^d)^q,\ g=K_1(g_0)\cdots K_q(g_{q-1})\right\}.$$
Notice in particular that any $g$ in ${\cal B}_q$ is bounded and inherits the regularity properties of the kernels $K_j$. Then, for $g\in{\cal B}_q$, $x\in\R^d$ and $L\in\R$, let us denote 
$$H_q^g(x,L)=\int_{S(x')=L}K_{q+1}(x,x')g(x')\frac{\dbar x'}{|DS(x')|}.$$

\textbf{Assumption [${\cal H}$]}
\begin{enumerate}
\item[$(i)$] For any $q\geq 0$, the mapping $x\mapsto H_q^1(x,L_q)$ belongs to $L^2(\eta)$, that is 
$$\int\eta(dx)\left(\int_{S(x')=L_q}K_{q+1}(x,x')\frac{\dbar x'}{|DS(x')|}\right)^2<\infty.$$
\item[$(ii)$] For any $q>0$, for any $g\in{\cal B}_q$, there exists $h\in L^2(\eta)$ such that for any $\varepsilon>0$, there exists $\delta>0$ such that for any $L\in[L_q-\delta,L_q+\delta]$ and for almost every $x\in\R^d$, 
$$\left|H_q^g(x,L)-H_q^g(x,L_q)\right|\leq \varepsilon h(x).$$
\end{enumerate}
\medskip

We will comment on this assumption in Section \ref{assumption}. In particular, we will see that it is not restrictive and is verified by most models of interest, for example when the level sets $\{S(x')=L\}$ have finite Hausdorff measure or when $\eta$ and the kernels $K_q$ have light tails.\medskip

If the kernels $K_j$ are based on Metropolis-Hastings algorithm as detailed in Section \ref{mh}, then one has to consider the set of functions
$${\cal B}_q^a=\left\{g:\R^d\to\R,\ \exists (g_0\dots,g_{q-1})\in\Ba(\R^d)^q,\ g=\ck^a_1(g_0)\cdots \ck^a_q(g_{q-1})\right\},$$
as well as the mapping
$$H_q^{g,a}(x,L)=\int_{S(x')=L}\ck^a_{q+1}(x,x')g(x')\frac{\dbar x'}{|DS(x')|},$$
and assumption [${\cal H}^a$] defined as follows.\medskip

\textbf{Assumption [${\cal H}^a$]}
\begin{enumerate}
\item[$(i)$] For any $q\geq 0$, the mapping $x\mapsto H_q^{1,a}(x,L_q)$ belongs to $L^2(\eta)$, that is 
$$\int\eta(dx)\left(\int_{S(x')=L_q}\ck^a_{q+1}(x,x')\frac{\dbar x'}{|DS(x')|}\right)^2<\infty,$$
and $\eta(H_q^{1,a}(.,L_q))>0$.
\item[$(ii)$] For any $q>0$, for any $g\in{\cal B}_q^a$, there exists $h\in L^2(\eta)$ such that for any $\varepsilon>0$, there exists $\delta>0$ such that for any $L\in[L_q-\delta,L_q+\delta]$ and for almost every $x\in\R^d$, 
$$\left|H_q^{g,a}(x,L)-H_q^{g,a}(x,L_q)\right|\leq \varepsilon h(x).$$
\end{enumerate}
\medskip

The main result of this paper is the following central limit type theorem.
  
\begin{theo}\label{theorem1}
Under Assumption $[{\cal H}]$ or $[{\cal H}^a]$, for any $f\in\Ba(\R^d)$ such that $f=f\times\un_{S(\cdot)\geq L^\star}$, we have
$$\sqrt{N}\left(\eta_n^N(f)-\eta_n(f)\right)\xrightarrow[N\to\infty]{\cal D}{\cal N}(0,\Gamma(f)),$$
with the variance functional
\begin{equation}\label{def-Gamma}
\Gamma(f):=\sum_{p=0}^n\e_p(\overline{Q}_{p,n}(f)^2-\e_n(f)^2).
\end{equation}
\end{theo}

Theorems \ref{kjsc} and \ref{theorem1} allow us to specify the fluctuations of the estimates $\hat{E}$, $\hat{P}$ and $\hat{C}$.

\begin{cor}\label{cor1}
Under the same assumptions as in Theorem \ref{theorem1}, we have:
\begin{enumerate}
\item[$(i)$] for the estimate of the expectation $E=\E[f(X)]=\E[f(X)\un_{S(X)\geq L^\star}]$,
$$\sqrt{N}\left(\hat{E}-E\right)\xrightarrow[N\to\infty]{{\cal D}}{\cal N}(0,\a^{2n}\Gamma(f)).$$
\item[$(ii)$] for the rare event probability $P=\P(Y\geq L^\star)$,
$$\sqrt{N}\left(\hat{P}-P\right)\xrightarrow[N\to\infty]{{\cal D}}{\cal N}(0,\a^{2n}\Gamma(\un_{S(\cdot)\geq L^\star})).$$
\item[$(iii)$] for the conditional expectation $C=\E[f(X)|S(X)\geq L^\star]$, still with $f=f\times\un_{S(\cdot)\geq L^\star}$, 
$$\sqrt{N}\left(\hat{C}-C\right)\xrightarrow[N\to\infty]{{\cal D}}{\cal N}(0,\Gamma(g)),$$
where 
$$g:=\frac{\un_{S(\cdot)\geq L^\star}}{r}\left(f-\frac{\e_{n}(f)}{r}\right).$$  
\end{enumerate}
\end{cor}

In the next section, we compare these results with the ones obtained for the fixed-levels version of Multilevel Splitting, which was initially proposed in \cite{alea06}. The analysis of this method in the specific context of the present article was done by some of the authors in \cite{cdfg}.

\subsection{Comparison with the fixed-levels method}\label{comparison}

In what follows, we return to the optimal Feynman-Kac particle approximation (fixed-levels method) that was presented in Section \ref{fx}.

\begin{theo}\label{delmo}
For any $f\in\Ba(\R^d)$, we have the almost sure convergences
$\lim_{N\rightarrow\infty}\check{\g}_n^N(f)=\g_n(f)$, and $\lim_{N\rightarrow\infty}\check{\e}_n^N(f)=\eta_n(f)$, as well as the 
convergences in distribution
\begin{eqnarray*}
\sqrt{N}\left(\check{\g}_n^N(f)-\g_n(f)\right)&\xrightarrow[N\to\infty]{\cal D}&{\cal N}(0,\a^{2n}\Gamma(f))\\
\sqrt{N}\left(\check{\e}_n^N(f)-\e_n(f)\right)&\xrightarrow[N\to\infty]{\cal D}&{\cal N}(0,\Gamma(f-\e_n(f)))
\end{eqnarray*}
with the variance functional $\Gamma$ defined in (\ref{def-Gamma}).
\end{theo} 

For the proof of this theorem, we report the interested reader to Propositions 9.4.1 and 9.4.2 in \cite{delmoral04a}. Just note that straightforward computations give that
\begin{eqnarray*}
\a^{2n}\Gamma(f)
&=&\a^{2n}\sum_{p=0}^n\e_p(\overline{Q}_{p,n}(f)^2-\e_n(f)^2)\\ \\
&=&\sum_{p=0}^n \gamma_p(1)^2 \eta_p((Q_{p,n}(f)-\eta_p(Q_{p,n}(f)))^2),
\end{eqnarray*}
which is exactly the variance given in Propositions 9.4.1  in \cite{delmoral04a} for the case of a multinomial resampling. In this paper we prefer using the first expression because it is how it will appear naturally in the proofs. In Theorem \ref{theorem1}, the asymptotic variance has the same form as the one for unnormalized measures in Theorem \ref{delmo} . Actually, in the adaptive case, we have $\gamma_n^N=\alpha^n \e_n^N$, so that both give the same asymptotic variance, up to a deterministic multiplicative constant.\medskip  

Note also that $\Gamma(f-\e_n(f))$ can be rewritten as 
$$
\Gamma(f-\e_n(f))=\sum_{p=0}^n\e_p(\overline{Q}_{p,n}(f-\e_n(f))^2),
$$
which is exactly formula (9.13) in \cite{delmoral04a}.\medskip  

In the normalized case, it may seem that Theorems \ref{theorem1} and \ref{delmo} give different asymptotic variances, and they do in all generality. But we need to carefully consider what we do at the last step and how we construct our estimates. The specificity of the last step is that the level is always $L^\star$ and thus is {\bf not} random.\medskip   


For the normalized measures, the asymptotic variances clearly coincide for functions $f$ such that $\eta_n(f)=0$. As we see in the proof of Corollary \ref{cor1}, we can write 
$$
\hat C  - C= \frac{\eta_n(\un_{S(\cdot)\geq L^\star})}{\eta_n^N(\un_{S(\cdot)\geq L^\star})}\times\left(\eta_n^N(g)-\eta_n(g)\right),
$$
where the prefactor $\frac{\eta_n(\un_{S(\cdot)\geq L^\star})}{\eta_n^N(\un_{S(\cdot)\geq L^\star})}$ converge to 1 in probability, and 
$$g=\frac{\un_{S(\cdot)\geq L^\star}}
{\eta_n(\un_{S(\cdot)\geq L^\star})}\left(f-\frac{\eta_n(f)}{\eta_n(\un_{S(\cdot)\geq L^\star})}\right)$$
 is such that $\eta_n(g)=0$. As the same trick can be done for $\check{C}-C$ (non adaptive case), we have the same asymptotic variance, because $g$ is centered for $\eta_n$.

\medskip
The next corollary, which is a direct consequence of Corollary~\ref{cor1} and the above discussion,  constitutes the main message of the present article.

\begin{cor}\label{comp}
Under Assumption $[{\cal H}]$ or $[{\cal H}^a]$, for any $f\in\Ba(\R^d)$ such that $f=f\times\un_{S(\cdot)\geq L^\star}$, the estimates $\hat{E}$ and $\check{E}$ have the same asymptotic variances.
The same result holds for the estimates $\hat{P}$ and $\check{P}$ of the probability $P$, and for the estimates $\hat{C}$ and $\check{C}$ of the conditional expectation $C$.
\end{cor}


Interestingly, as detailed in Proposition 3 of \cite{cdfg}, there exists another expression for the asymptotic variance of the estimator $\check{P}$. By Corollary \ref{comp}, this expression holds for the estimator $\hat{P}$ as well. We recall it hereafter for the sake of completeness.

\begin{cor}\label{zijd}
Under Assumption $[{\cal H}]$ or $[{\cal H}^a]$, we have
$$\sqrt N\ \frac{\hat P-P}{P}\xrightarrow[N\to+\infty]{{\cal D}}{\cal N}(0,\sigma^2)\quad\mbox{and}\quad\sqrt N\ \frac{\check P-P}{P}\xrightarrow[N\to+\infty]{{\cal D}}{\cal N}(0,\sigma^2),$$
where $\sigma^2=\frac{\a^{2n}}{P^2}\Gamma(\un_{S(\cdot)\geq L})$ admits the alternative expression
\begin{align}
\sigma^2=&\ (n-1)\times\frac{1-\a}{\a}+\frac{1-r}{r}\nonumber\\
&+\frac{1}{\a}\sum_{p=0}^{n-2}\E\left[\left.\left(\frac{\P(S(X_n)\geq L^\star|X_{p+1})}{r\times\alpha^{n-(p+1)}}-1\right)^2\right|S(X_{p})\geq L_{p}\right]\nonumber\\
&+\frac{1}{r}\times\E\left[\left.\left(\frac{\P(S(X_n)\geq L^\star|X_n)}{r}-1\right)^2\right|S(X_{n-1})\geq L_{n-1}\right].\label{rv}
\end{align}
\end{cor}

This expression emphasizes that, when using Multilevel Splitting, the relative variance $\sigma^2$ is always lower bounded by an incompressible variance term, namely that
$$\sigma^2\geq (n-1)\times\frac{1-\a}{\a}+\frac{1-r}{r}.$$ The additive terms in (\ref{rv}) depend on the mixing properties of the transition kernels $M_p$. In particular, if at each step we have an ``ideal'' kernel, meaning that, knowing that $S(X_p)\geq L_p$, $X_{p+1}$ is independent of $X_p$, then these additive terms vanish. This is the so-called ``idealized'' version of Adaptive Multilevel Splitting, studied for example in \cite{cdfg,ghm,blr14,s14}.
\medskip

 Finally, let us mention that our results also apply directly to the case of general multilevel splitting. Specifically, let us consider a fixed and known final level $L^\star$ and a sequence of prescribed success probabilities $(\alpha_p)_{p\geq 0}$ corresponding to the (unknown) sequence of levels $-\infty = L_{-1}<L_0 < \dots < L_{n-1}<L^\star<L_n$, with 
$$\alpha_p=\P(S(X) > L_{p} | S(X) > L_{p-1})\ \Longrightarrow\ \P(S(X) > L^\star) =r \prod_{p=0}^{n-1}\alpha_p,$$
with 
$$r=\P(S(X) > L^\star | S(X) > L_{n-1})\in(\alpha_n,1).$$ 
Then, a quick inspection of the proofs ensures that the Adaptive Multilevel Splitting algorithm with a sequence of adaptive levels $L_0^N < \dots < L_{n-1}^N$ will have the same asymptotic variance as the Multilevel Splitting algorithm with the levels $L_0 < \dots < L_{n-1}$. Compared to its fixed-levels counterpart, the cost of the adaptive version is just a higher complexity by a factor $\log N$, due to the quicksort of the sample at each step.

\section{Discussion on Assumption [${\cal H}$]}\label{assumption}
In this section we return to Assumption [${\cal H}$], and show that it is satisfied on several models of interest. For the sake of simplicity, we focus our attention on [${\cal H}$] and we will not comment on [${\cal H}^a$], but the following arguments may be repeated {\it mutatis mutandis} by replacing $K_q$ with $\ck_q^a$.

\subsection{An alternative formulation of [${\cal H}$]}
 
First we explain how Assumption [${\cal H}$]$(ii)$ can be verified via a condition on the kernels $K_q$. Specifically, we get an assumption which is easier to check than [${\cal H}$]$(ii)$, while only a bit more restrictive. The proof is given in Section \ref{preuveH}.   

\begin{pro}
\label{nec.cond}
Assumption $[{\cal H}](ii)$ is fulfilled if there exists a function $h$ in $L^2(\eta)$ and a real number $\delta>0$, such that for any $q>0$, for all $L\in(L_q-\delta,L_q+\delta)$ and for almost every $x\in\R^d$,
\begin{align}
&\int_{S(x')=L} \left| \div_{x'}\left(K_{q+1}(x,x')\ \frac{DS(x')}{|DS(x')|^2}\right) \right| \frac{\dbar x'}{|DS(x')|}\nonumber\\
&\ +\int_{S(x')=L} \left[\sum_{m=1}^q  \sum_{j=1}^d \int\left|\frac{\partial}{\partial x_j'}K_{m}(x',x'')\right|dx''\right]K_{q+1}(x,x')\frac{\dbar x'}{|DS(x')|^2}\leq h(x).\label{eq.H2}
\end{align}
\end{pro}

One may think at first sight that this condition is not much easier to handle than $[{\cal H}](ii)$ but, considering our framework, we stress the fact that it is much more natural since it involves only the measure $\eta$, the score function $S$ and the transition kernels $K_q$, at the cost of only a very slight restriction. We can also remark that when $K_m=K$ for all $m$, then the sum in $m$ in equation~ (\ref{eq.H2}) can be omitted.

\subsection{Examples}

This section exhibits two typical situations where Assumption [${\cal H}$] is satisfied. In order to verify $[{\cal H}](ii)$, we will make use of condition (\ref{eq.H2}) given in Proposition \ref{nec.cond}.

\subsubsection{The compact case} 

If the levels sets $S^{-1}(L)$ are compact, then under mild regularity conditions on $S$ and the kernels $K_q$, it is not difficult to see that Assumption [${\cal H}$] is satisfied. The remainder of this section details some sufficient conditions.\medskip

Let us assume that, for any real number $L$, the set $\{S(x)\leq L\}$ is bounded, with Hausdorff measure bounded by $C_L$, and that for all $x\in S^{-1}(L)$, we have $|DS(x)|\geq c_L>0$. Besides, assume that for all $q$, $K_q\leq C$. From these assumptions, it is clear that 
$$
\int\eta(dx)\left(\int_{S(x')=L_q}K_{q+1}(x,x')\frac{\dbar x'}{|DS(x')|}\right)^2\leq \left(\frac{C_q\times C}{c_q}\right)^2,
$$
and condition [${\cal H}$]$(i)$ is fulfilled.
\medskip

Now we consider [${\cal H}$]$(ii)$, and use condition (\ref{eq.H2}) of Proposition \ref{nec.cond}. Denote respectively by $C_{q,\delta}$ and $c_{q,\delta}$ the upper and lower bounds of $C_L$ and $c_L$ for $L\in (L_q-\delta,L_q+\delta)$. We reformulate the first term in (\ref{eq.H2}) as the integral on the level set $\{S(x')=L\}$ of the quantity
\begin{eqnarray*}
\lefteqn{
\div_{x'}\left[\frac{K_{q+1}(x,x')}{|DS(x')|^2} DS(x')\right] }\\
&=&\frac{K_{q+1}(x,x')}{|DS(x')|^2}\ \Delta S(x') + \frac{1}{|DS(x')|^2}\  D_{x'}K_{q+1}(x,x')DS(x')^T\\
&&- 2\ \frac{K_{q+1}(x,x')}{|DS(x')|^4}\ DS(x')H_S(x')DS(x')^T,
\end{eqnarray*}
where $\Delta S$ and $H_S$ are respectively the Laplacian and the Hessian of $S$. From this we see that if $K_q$ has bounded first derivatives (in the second variable), and if $S$ is two times continuously differentiable, then this term is bounded by a constant $M$ and, returning to (\ref{eq.H2}), we get for all $L\in (L_q-\delta,L_q+\delta)$, 
$$\int_{S(x')=L} \left| \div_{x'}\left(K_{q+1}(x,x')\ \frac{DS(x')}{|DS(x')|^2}\right) \right| \frac{\dbar x'}{|DS(x')|}\leq M\times\frac{C_{q,\delta}}{c_{q,\delta}}.$$
The second term in (\ref{eq.H2}), namely
$$\int_{S(x')=L} \left[\sum_{m=1}^q  \sum_{j=1}^d \int\left|\frac{\partial}{\partial x_j'}K_{m}(x',x'')\right|dx''\right]K_{q+1}(x,x')\frac{\dbar x'}{|DS(x')|^2},$$
is slightly more challenging because of the inner integral on the whole space. An obvious sufficient condition is that the kernels $K_m$ have bounded first derivatives in the first variable, say by $M$, and that their supports have uniformly bounded Lebesgue measures, say by $\rho$. Then we have, for all $L\in (L_q-\delta,L_q+\delta)$,  
$$\int_{S(x')=L}\left[\sum_{m=1}^q\sum_{j=1}^d\int\left|\frac{\partial}{\partial x_j'}K_{m}(x',x'')\right|dx''\right]K_{q+1}(x,x')\frac{\dbar x'}{|DS(x')|^2}\leq qd\rho MC\frac{C_{q,\delta}}{c_{q,\delta}^2}$$
and condition (\ref{eq.H2}) of Proposition \ref{nec.cond} is satisfied in this so-called compact case (compact level sets for $S$ plus compact supports for the transition kernels).

\subsubsection{The Gaussian case}

Outside this compact framework, there are of course other situations where Assumption [${\cal H}$] is satisfied. Indeed, in many cases, $K_q$ and $\eta$ have exponential decay at infinity (i.e. light tails). If $K_q$ has the form
$$K_q(x,x')\propto\exp(-V_q(x-x')),$$ 
with $V_q$ equivalent to a polynomial at infinity, then basically its derivatives with respect to $x'$ take the form $P_q(x,x')K_q(x,x')$, where $P_q$ itself is bounded by another polynomial at infinity. Then, roughly speaking, its integral is a moment of the density $K_q(x,\cdot)$, which typically will be bounded by another polynomial in $x$. This polynomial in $x$ will in turn be integrable by $\eta$ if $\eta$ has light tails. The upcoming example is going to make this more precise in the Gaussian case. 
\medskip 

Specifically, we will detail the computations on the zero-bit watermarking example of \cite[section~5.1]{cdfg}. In this case, the score function is defined for any $x\in\R^d$ by $S(x)=x_1/|x|$, and $\eta$ is the standard Gaussian distribution on $\R^d$. Thus it is readily seen that, for any $\sigma>0$, the transition kernel $K$ defined by
$$K(x,x')=\frac{1+\s^2}{2\pi\s^2}\exp\left(-\frac{1+\s^2}{2\s^2}\left|x'-\frac{x}{\sqrt{1+\s^2}}\right|^2\right)$$
is $\eta$-reversible. We explain in Section \ref{preuvegauss} that in this situation, Assumption [${\cal H}$] is satisfied.

\section{Proofs}\label{azp}

\subsection{Some preliminary notations}\label{prel}

We let ${\cal F}_{-1}^N:=\{\emptyset,\Omega\}$ be the trivial sigma-field and, for $q\geq 0$, we denote by ${\cal F}_{q}^N$ 
the sigma-field generated 
$${\cal F}_{q}^N:={\cal F}_{q-1}^N\vee\sigma\left((X_{q}^1,U_{q}^1),\dots,(X_{q}^N,U_{q}^N)\right).$$
Also, according to the definition of the empirical quantile given by (\ref{defquant}), we set
$${\cal G}_{-1}^N:=\sigma\left(L_{0}^N\right)=\sigma\left(X_{0}^{(\lfloor N(1-\alpha)\rfloor)},U_{0}^{(\lfloor N(1-\alpha)\rfloor)}\right),$$ 
and, for $q\geq 0$,
$${\cal G}_{q}^N:={\cal F}_{q}^N\vee\sigma\left(L_{q+1}^N\right)={\cal F}_{q}^N\vee\sigma\left(X_{q+1}^{(\lfloor N(1-\alpha)\rfloor)},U_{q+1}^{(\lfloor N(1-\alpha)\rfloor)}\right).$$
Then, given ${\cal F}_{q-1}^N$, 
$$\eta^N_q:=\frac{1}{N}\sum_{1\leq i\leq N}\delta_{X^i_q}$$
is the empirical measure associated with $N$ conditionally independent random vectors with common distribution
$$
\Phi_q(\eta^N_{q-1}):=\frac{1}{\lceil N\alpha\rceil}\sum_{i:~S(X_{q-1}^i)>L^N_{q-1}}M_{q,L^N_{q-1}}(X_{q-1}^i,\cdot).
$$
Next, given ${\cal G}_{q-1}^N$ and adapting for instance Theorem 2.1 in~\cite{bcg-2012} to our context, it can be shown that the subsample of the vectors $X^i_q$ above $L_q^N$ are conditionally independent random vectors denoted by $(\widetilde{X}^i_q)_{1\leq i\leq\lceil N\alpha\rceil}$ and with common distribution
\begin{equation}\label{decor-refbis}
\mbox{\rm Law}\left(\left.\left(\widetilde{X}^1_q,\ldots,\widetilde{X}^{\lceil N\alpha\rceil}_q\right)~\right|~{\cal G}_{q-1}^N\right)=\Psi_{G_{\eta^N_{q}}}\left(\Phi_q(\eta^N_{q-1})\right)^{\otimes \lceil N\alpha\rceil},
\end{equation}
where, if $L_q^N=(L,u)$, we have for any $x\in\R^d$  
$$G_{\eta^N_{q}}(x)=\un_{S(x)>L}+(1-u)\un_{S(x)=L},$$
and, accordingly,
\begin{equation}\label{decor-reftetra}
\Psi_{G_{\eta^N_{q}}}\left(\Phi_q(\eta^N_{q-1})\right)=\frac{G_{\eta^N_{q}}\times\Phi_q(\eta^N_{q-1})}{\Phi_q(\eta^N_{q-1})(G_{\eta^N_{q}})}.
\end{equation}
In summary,  we have that
\begin{equation}\label{decor-ref}
\widetilde{\eta}^N_q:=\Psi_{G_{\eta^N_q}}\left(\eta^N_q\right)=\frac{1}{\lceil N\alpha\rceil}~\sum_{i=1}^{\lceil N\alpha\rceil}\delta_{\widetilde{X}^i_q},
\end{equation}
and
\begin{equation}\label{decor-refter}
\Phi_{q+1}\left(\eta^N_q\right)=\Psi_{G_{\eta^N_{q}}}\left(\eta^N_{q}\right)M_{q+1,L^N_{q-1}}=\widetilde{\eta}^N_qM_{q+1,L_q^N}.
\end{equation} 
Let us also define
\begin{equation}\label{yayaya}
\Pi_q(\eta^N_{q-1}):=\frac{\lceil N\alpha\rceil}{N}\times \Phi_q(\eta^N_{q-1})=\frac{\lceil N\alpha\rceil}{N}\times\Psi_{G_{\eta^N_{q-1}}}\left(\eta^N_{q-1}\right)M_{q,L^N_{q-1}}.
\end{equation}
Alternatively, if $\nu$ is absolutely continuous, we define the operator $\Pi_q$ (see also (\ref{lknc})) as 
\begin{equation}\label{yayayaya}
\Pi_q(\nu):=\alpha\ \Phi_{q}(\nu)=\alpha\ \Psi_{G_{\nu}}(\nu)M_{q,\nu}.
\end{equation} 
Besides, for any $q<p$ and $\mu=\nu$ or $\mu=\eta_q^N$, we set
\begin{equation}\label{yoyo}
\Pi_{q,p}(\mu)=\Pi_{q+1}(\mu)Q_{q+1,p}
\quad\mbox{\rm and}\quad
\Phi_{q,p}(\mu)=\frac{\Phi_{q+1}(\mu)Q_{q+1,p}}{(\Phi_{q+1}(\mu)Q_{q+1,p})(1)},
\end{equation}
with the conventions that $\Pi_{q,p}=I_d=\Phi_{q,p}$ whenever $q\geq p$. This yields
\begin{equation}\label{pi}
\e_p=\a^{q-p}\times \Pi_{q,p}(\e_q)=\Phi_{q,p}(\eta_q).
\end{equation}
Hence, for any  $f\in\Ba(\R^d)$, we have
\begin{equation}\nonumber
\begin{array}{l}
\Pi_{q,p}(\mu)
=\mu Q_{q,p,\mu}\quad\mbox{\rm and}\quad \Phi_{q,p}(\mu)(f)={\Pi_{q,p}(\mu)(f)}/{\Pi_{q,p}(\mu)(1)},
\end{array}
\end{equation}
with the collection of integral operators $Q_{q,p,\mu}$ defined by
\begin{eqnarray}\nonumber
Q_{q,p,\mu}&:=&
Q_{q+1,\mu}Q_{q+2}\dots Q_{p}=Q_{q+1,\mu}Q_{q+1,p}.
\end{eqnarray}
In addition, using (\ref{decor-ref}), we prove 
\begin{equation}\label{io}
\Pi_{q+1}(\eta_{q}^N)=\eta^N_qQ_{q+1,\eta^N_q}=\eta^N_q(G_{\eta^N_q})~\Phi_{q+1}(\eta^N_q)=
\frac{\lceil N\alpha\rceil}{N}~~\widetilde{\eta}^N_qM_{q+1,L_q^N},
\end{equation}
which implies that
\begin{equation}\label{bla}
\Pi_{q,p}(\eta^N_q)=\frac{\lceil N\alpha\rceil}{N}~~\widetilde{\eta}^N_q\widetilde{Q}_{q,p,\eta^N_q}
\end{equation}
whence, thanks to (\ref{decor-refbis}), 
\begin{equation}\label{ref-cond-G}
\E\left[\Pi_{q,p}(\eta^N_q)(f)~|~{\cal G}_{q-1}^N\right]=\frac{\lceil N\alpha\rceil}{N}~\Psi_{G_{\eta^N_q}}\left(\Phi_q(\eta^N_{q-1})\right)\widetilde{Q}_{q,p,\eta^N_q}(f),
\end{equation}
with the collection of integral operators
\begin{equation}\nonumber
\widetilde{Q}_{q,p,\mu}:=M_{q+1,\mu}Q_{q+1,p}.
\end{equation}
Note that by construction, we have
$$
\widetilde{Q}_{q,p}=M_{q+1}Q_{q+1,p}=:\widetilde{Q}_{q,p}
\quad\mbox{\rm and}\quad
Q_{q,p,\eta_q}=Q_{q,p}.
$$
We also observe that, according to (\ref{yayaya}),
\begin{equation}\label{Tq-Phi0}
\E\left[\eta_q^N(f)~|~{\cal F}_{q-1}^N\right]=\Phi_{q }\left(\eta_{q-1}^N\right)(f)=\frac{N}{\lceil N\alpha\rceil}\Pi_{q}(\eta_{q-1}^N)(f),
\end{equation}
or, said differently,
\begin{equation}\label{Tq-Phi}
\alpha^{-1}~\Pi_{q}(\eta_{q-1}^N)=
\rho_N~\Phi_{q }\left(\eta_{q-1}^N\right)\quad\mbox{\rm with}\quad
\rho_N:=\frac{\lceil N\alpha\rceil}{N\alpha}
\end{equation}
and
\begin{equation}\label{Tq-Phi2}
\alpha^{-1}~\Pi_{q-1,p}\left(\eta_{q-1}^N\right)=\alpha^{-1}~\Pi_{q,p}\left(\Pi_{q}(\eta_{q-1}^N)\right)=\rho_N~~\Pi_{q,p}\left(\Phi_{q }\left(\eta_{q-1}^N\right)\right).
\end{equation}
We note, once and for all, that 
$$0\leq \rho_N-1<\frac{1}{N\alpha}.$$
Finally, we consider the ${\cal G}_{q-1}^N$ measurable random variable $\epsilon_{q}^N$ defined by
\begin{equation}\label{yo}
\epsilon_{q}^N=1-\rho_N\Phi_{q }\left(\eta_{q-1}^N\right)(G_{\eta_{q}^N})/\alpha\ \Longleftrightarrow\ \Phi_{q }\left(\eta_{q-1}^N\right)(G_{\eta_{q}^N})=\rho_N^{-1}\alpha(1-\epsilon_{q}^N).
\end{equation}

\subsection{Proof of Theorem \ref{kjsc}}
We will prove the almost sure convergences, and explain at the end how to get the convergence in probability. We proceed by induction with respect to the time parameter $p$, as is done for example in \cite{chopin2004,doucmoulines}.
\medskip
 
Denoting $X_0^1,\dots,X_0^N$ an i.i.d. sample with common law $\eta=\eta_0$, the strong law of large numbers tells us that, by definition of $\eta_0^N$ and $\eta_0$, for any $f\in L^2(\eta)$, we have
$$\e_0^N(f)=\frac{1}{N}\sum_{i=1}^{N}f(X_0^i)\xrightarrow[N\to\infty]{a.s.}\e(f)=\e_0(f).$$
Then, since the cdf $F_Y$ is one-to-one and $L_0=F_Y^{-1}(1-\a)$, the theory of order statistics ensures that  
$$L_0^N=L_{\eta^N_0}\xrightarrow[N\to\infty]{a.s.}L_{\eta_0}=L_0.$$
Next, let us assume that the property is satisfied for $p\geq 0$ and recall that ${\cal F}_p^N$ is the sigma-field generated by the random couples $(X_p^i,U_p^i)$ for $i=1,\dots,N$. We begin with the following decomposition
\begin{align}
&\left|\e_{p+1}^N(f)-\e_{p+1}(f)\right|\nonumber\\
&\quad\leq\left|\e_{p+1}^N(f)-\E\left[\left.\e_{p+1}^N(f)\right|{\cal F}_p^N\right]\right|+\left|\E\left[\left.\e_{p+1}^N(f)\right|{\cal F}_p^N\right]-\e_{p+1}(f)\right|.\label{oihx}
\end{align}
Concerning the second term, (\ref{Tq-Phi0}) implies
$$\E\left[\left.\e_{p+1}^N(f)\right|{\cal F}_p^N\right]=\Phi_{p+1}(\eta^N_{p})(f)=\frac{N}{\lceil N\alpha\rceil}~ \Pi_{p+1}(\eta^N_{p})(f),$$
and by Proposition \ref{asc} page \pageref{asc}, the induction assumption and (\ref{pi}), we get
$$\frac{N}{\lceil N\alpha\rceil}~\Pi_{p+1}(\eta^N_{p})(f)\xrightarrow[N\to\infty]{a.s.}\frac{1}{\alpha}~\Pi_{p+1}(\eta_{p})(f)=\eta_{p+1}(f).$$
Hence the second term of (\ref{oihx}) goes almost surely to 0. For the first term of (\ref{oihx}), recall that given ${\cal F}_p^N$, the random variables $f(X_{p+1}^1),\dots,f(X_{p+1}^N)$ are i.i.d. with mean $\E\left[\left.\e_{p+1}^N(f)\right|{\cal F}_p^N\right]$. Hence, for any $\varepsilon>0$, Hoeffding's inequality gives
\begin{align}\label{hoeffding}
&\P\left(\left.\left|\e_{p+1}^N(f)-\E\left[\left.\e_{p+1}^N(f)\right|{\cal F}_p^N\right]\right|>\varepsilon\right|{\cal F}_p^N\right)\leq 2\exp\left\{-\frac{N\varepsilon^2}{2\|f\| ^2}\right\},
\end{align}
Since this upper-bound is deterministic, this amounts to say that
\begin{align}\nonumber
&\P\left(\left|\e_{p+1}^N(f)-\E\left[\left.\e_{p+1}^N(f)\right|{\cal F}_p^N\right]\right|>\varepsilon\right)\leq 2\exp\left\{-\frac{N\varepsilon^2}{2\|f\| ^2}\right\}.
\end{align}
Consequently, the choice $\varepsilon_N=N^{-1/4}$  and Borel-Cantelli Lemma show that the first term of (\ref{oihx}) goes almost surely to 0 as well.\medskip

It remains to show the convergence of $L_{p+1}^N$ to $L_{p+1}$. To achieve this aim, let us denote $F_{p+1}$ the following cdf
$$F_{p+1}(y)=\P(S(X)\leq y~|~S(X)\geq L_{p}).$$
In this respect, by definition, we have $F_{p+1}(L_{p+1})=1-\a$. This being done, one has just to mimic the reasoning of the proof of point $(i)$ in Proposition \ref{asc} to obtain the desired result. 
\medskip

To get the convergences in probability for functions $f\in L^2(\eta)$, the same arguments apply to the second term of (\ref{oihx}). About the first one, one may just replace Hoeffding's inequality with Chebyshev's inequality in (\ref{hoeffding}) to obtain
\begin{align}\label{chebyshev}
&\P\left(\left.\left|\e_{p+1}^N(f)-\E\left[\left.\e_{p+1}^N(f)\right|{\cal F}_p^N\right]\right|>\varepsilon\right|{\cal F}_p^N\right)\leq \frac{\sigma_N^2}{\varepsilon^2},
\end{align}
where 
$$\sigma_N^2=\E\left[\left.\left(\e_{p+1}^N(f)-\E\left[\left.\e_{p+1}^N(f)\right|{\cal F}_p^N\right]\right)^2\right|{\cal F}_p^N\right].$$ 
Given ${\cal F}_p^N$, the random variables $X_{p+1}^1,\dots,X_{p+1}^N$ are i.i.d. with law $\Phi_{p+1}(\eta_p^N)$, so 
$$\sigma_N^2=\frac{1}{N}\left\{\Phi_{p+1}(\eta_p^N)(f^2)-\Phi_{p+1}(\eta_p^N)(f)^2\right\}.$$ 
Obviously, by (\ref{Tq-Phi0}), the induction assumption and Proposition \ref{asc},
$$
N\sigma_N^2\leq \Phi_{p+1}(\eta_p^N)(f^2)= \frac{N}{\lceil N\alpha\rceil}~\Pi_{p+1}(\eta^N_{p})(f^2)\xrightarrow[N\to\infty]{\P}\eta_{p+1}(f^2). 
$$
This proves that, for any $\varepsilon>0$, 
$$\P\left(\left.\left|\e_{p+1}^N(f)-\E\left[\left.\e_{p+1}^N(f)\right|{\cal F}_p^N\right]\right|>\varepsilon\right|{\cal F}_p^N\right)\xrightarrow[N\to\infty]{\P}0,$$
and Lebesgue's dominated convergence ensures that
$$\e_{p+1}^N(f)-\E\left[\left.\e_{p+1}^N(f)\right|{\cal F}_p^N\right]\xrightarrow[N\to\infty]{\P}0.$$
This concludes the proof of Theorem \ref{kjsc}.
\hfill  $\blacksquare$ 

\subsection{Proof of Theorem \ref{theorem1}}\label{pozd}

We use the symbols $\V(\point)$ and $\V(\point~|~{\cal G}_{q}^N)$ to denote respectively the variance and the conditional variance operators. We start the analysis with a decomposition which is equivalent to the one given  (for example) in \cite{delmoral04a} page 216. Specifically, for any $p\geq 0$, we have the standard following telescoping sum
$$
\eta_p^N-\eta_p=\sum_{q=0}^p\alpha^{q-p}\left\{\Pi_{q,p}(\eta_q^N)-\alpha^{-1}~\Pi_{q-1,p}(\eta_{q-1}^N)\right\},\\
$$
with the conventions $\eta_{-1}^N=\eta_0=\eta$ and $\Pi_{0}=\alpha I_d$.  By (\ref{Tq-Phi2}), this implies that
\begin{equation}\label{pozdk}
[\eta_p^N-\eta_p](f)\\
={\cal M}_p^N + {\cal R}_p^N,
\end{equation}
where 
\begin{eqnarray}\label{eqM1}
\lefteqn{{\cal M}_p^N=
\sum_{q=0}^p\alpha^{q-p}\left\{\Pi_{q,p}(\eta_q^N)(f)-\E\left[\left.\Pi_{q,p}(\eta_q^N)(f)\right|{\cal G}_{q-1}^N\right]\right\}}\\
&&\nonumber+\sum_{q=0}^p\alpha^{q-p}\left\{
\epsilon_{q}^N~\E\left[\left.\Pi_{q,p}(\eta_q^N)(f)\right|{\cal G}_{q-1}^N\right]-
\E\left[\left.\epsilon_{q}^N~\E\left[\left.\Pi_{q,p}(\eta_q^N)(f)\right|{\cal G}_{q-1}^N\right]
\right|{\cal F}_{q-1}^N\right]
\right\}
\end{eqnarray}
is a martingale that will be discussed below, and
\begin{eqnarray}\label{eqR1}
 {\cal R}_p^N &=&\sum_{q=0}^p\alpha^{q-p}
 \E\left[\left.\epsilon_{q}^N~\E\left[\left.\Pi_{q,p}(\eta_q^N)(f)\right|{\cal G}_{q-1}^N\right]\right|{\cal F}_{q-1}^N\right]
\\\nonumber
&&+\sum_{q=0}^p\alpha^{q-p}\left\{(1-\epsilon_{q}^N)\E\left[\left.\Pi_{q,p}(\eta_q^N)(f)\right|{\cal G}_{q-1}^N\right]-\alpha^{-1}~\Pi_{q-1,p}(\eta_{q-1}^N)\right\},
\end{eqnarray}
is a rest that will be negligible. We recall that $\rho_N=\lceil N\alpha\rceil/N$ and  that $\epsilon_{q}^N$ was defined in equation (\ref{yo}) by
$$\epsilon_{q}^N=1-\frac{\rho_N}{\alpha}\times\Phi_{q }\left(\eta_{q-1}^N\right)(G_{\eta_{q}^N}).$$
The analysis of (\ref{pozdk}) is based on a series of technical results.

\begin{pro}\label{oisx}
For any $q\leq p$ and any $f\in\Ba(\R^d)$, we have
$$
N\a^{2(q-p)}~\V\left(\left.\Pi_{q,p}(\eta_q^N)(f)\right|{\cal G}_{q-1}^N\right)\xrightarrow[N\to\infty]{a.s.}\e_q(\overline{Q}_{q,p}(f)^2)-\a^{-1}~\e_p(f)^2.\nonumber
$$
\end{pro}

\begin{pro}\label{oisw}
For any $q\leq p$ and any $f\in\Ba(\R^d)$, we have
$$
\E\left[\left.\sqrt{N}~\epsilon_{q}^N~\E\left[\left.\Pi_{q,p}(\eta_q^N)(f)\right|{\cal G}_{q-1}^N\right]\right| {\cal F}_{q-1}^N\right]\xrightarrow[N\to\infty]{\P} 0,
$$
and 
$$
\V\left[\left.\sqrt{N}~\epsilon_{q}^N~\E\left[\left.\Pi_{q,p}(\eta_q^N)(f)\right|{\cal G}_{q-1}^N\right]\right| {\cal F}_{q-1}^N\right]
 \xrightarrow[N\to\infty]{\P} \frac{1-\alpha}{\alpha}\eta_p(f)^2.
$$
\end{pro}

\begin{pro}\label{VieillePuteAlbanaise}
Under Assumption $[{\cal H}]$, for any $q\leq p$ and any $f\in\Ba(\R^d)$ such that $f=f\times\un_{S(\cdot)\geq }L^\star$, we have
$$\sqrt{N}~\left(\alpha~(1-\epsilon_{q}^N)\E\left[\left.\Pi_{q,p}(\eta_q^N)(f)\right|{\cal G}_{q-1}^N\right]-\Pi_{q-1,p}(\eta_{q-1}^N)(f)\right)\xrightarrow[N\to\infty]{\P}0.$$
\end{pro}

The proofs of these propositions are detailed in Section \ref{technical}. Now we return to the proof of Theorem \ref{theorem1} by considering the decomposition (\ref{pozdk}). By Propositions \ref{oisw} and \ref{VieillePuteAlbanaise}, we have that 
$$\sqrt{N}~{\cal R}_p^N \xrightarrow[N\to\infty]{\P} 0.$$
By (\ref{bla}), (\ref{decor-ref})  and (\ref{Tq-Phi}), we may write ${\cal M}_p^{N}={\cal M}_p^{N,1}+{\cal M}_p^{N,2}$ with
\begin{align}
{\cal M}_p^{N,1}=&\sum_{q=0}^p\alpha^{q-p}\left(\Pi_{q,p}(\eta_q^N)(f)-\E\left[\left.\Pi_{q,p}(\eta_q^N)(f)\right|{\cal G}_{q-1}^N\right]\right)\nonumber\\
=&\frac{\lceil N\alpha\rceil}{N}~ \sum_{q=0}^p\alpha^{q-p}\left(\widetilde{\eta}^N_q\widetilde{Q}_{q,p,\eta^N_q}(f)-\E\left[\left.\widetilde{\eta}^N_q\widetilde{Q}_{q,p,\eta^N_q}(f)\right|{\cal G}_{q-1}^N\right]\right)\nonumber\\
=&\frac{\rho_N}{\lceil N\alpha\rceil}~ \sum_{q=0}^p\sum_{i=1}^{\lceil N\alpha\rceil}\alpha^{q-p+1}\left(\widetilde{Q}_{q,p,\eta^N_q}(f)(\widetilde{X}^i_q)-\E\left[\left.\widetilde{Q}_{q,p,\eta^N_q}(f)(\widetilde{X}^i_q)\right|{\cal G}_{q-1}^N\right]\right)\nonumber        
\end{align}
and
\begin{align}
&{\cal M}_p^{N,2}=\nonumber\\
&~\sum_{q=0}^p\alpha^{q-p}\left(\epsilon^N_q~\E\left[\left.\Pi_{q,p}(\eta_q^N)(f)\right|{\cal G}_{q-1}^N\right]-\E\left[\left.\epsilon_{q}^N~\E\left[\left.\Pi_{q,p}(\eta_q^N)(f)\right|{\cal G}_{q-1}^N\right]
\right|{\cal F}_{q-1}^N\right]\right).\nonumber
\end{align}
Now, remember the role of the auxiliary variables $(U_q^1,\dots,U_q^N)$ as mentioned on pages \pageref{ppu1} and \pageref{decor-refbis}, and consider the filtration ${\cal J}=({\cal J}_j)_{0\leq j\leq (p+1)(\lceil N\alpha\rceil+1)-1}$  constructed as
follows: for $q\in\{0,\dots,p\}$,
$${\cal J}^N_{q(\lceil N\alpha\rceil+1)}={\cal G}_{q-1}^N$$ 
and for $q\in\{0,\dots,p\}$ and $i\in\{1,\dots,\lceil N\alpha\rceil\}$,
$${\cal J}^N_{q(\lceil N\alpha\rceil+1)+i}={\cal G}_{q-1}^N\vee\sigma((\tilde
 X_q^{1},\tilde U_q^1),\dots,(\tilde  X_q^{i},\tilde U_q^i)) \vee \Sigma_q^N$$
 with
 $$ \Sigma_q^N = \sigma((X_q^{j}, U_q^j), \mbox{ for the indices $j$ such that $S(X^j_q)<L_q^N$ }).$$
In particular, note that 
$${\cal J}_{q(\lceil N\alpha\rceil+1)+\lceil N\alpha\rceil}^N={\cal J}^N_{(q+1)(\lceil N\alpha\rceil+1)-1}={\cal F}_q^N.$$ 
Let us define the sequence of random variables $(Z_j^N)_{0\leq j\leq (p+1)(\lceil N\alpha\rceil+1)-1}$ where the term of rank $q(\lceil N\alpha\rceil+1)$ is
$$
\begin{array}{l}
{Z^N_{q(\lceil N\alpha\rceil+1)}}\\
\\
=\alpha^{q-p}\left\{\epsilon^N_q~\E\left[\left.\Pi_{q,p}(\eta_q^N)(f)\right|{\cal G}_{q-1}^N\right]-\,\E\left[\left.\epsilon^N_q~\E\left[\left.\Pi_{q,p}(\eta_q^N)(f)\right|{\cal G}_{q-1}^N\right]\right|{\cal F}_{q-1}^N\right]\right\}\\ \\
=\alpha^{q-p}\left\{\epsilon^N_q~\E\left[\left.\Pi_{q,p}(\eta_q^N)(f)\right|  {\cal J}^N_{q(\lceil N\alpha\rceil+1)}
 \right]-\right.\\ \\
 \quad\quad\quad\quad\left.\E\left[\left.\epsilon^N_q~\E\left[\left.\Pi_{q,p}(\eta_q^N)(f)\right| {\cal J}^N_{q(\lceil N\alpha\rceil+1)} 
 \right]\right| {\cal J}^N_{q(\lceil N\alpha\rceil+1)-1} \right]\right\},
\end{array} 
$$
while the term of rank $ q(\lceil N\alpha\rceil+1)+i$, with $1\leq i\leq \lceil N\a\rceil$, is
$$
\begin{array}{l}
Z^N_{q(\lceil N\alpha\rceil+1)+i}\\
\\
=\frac{\rho_N}{\lceil N\alpha\rceil}~
\a^{q-p+1}\times\left\{
~\widetilde{Q}_{q,p,\eta^N_q}(f)(\tilde X_q^{i})-
\E\left[\left.\widetilde{Q}_{q,p,\eta^N_q}(f)(\tilde X_q^{i})\right|{\cal
    J}_{q(\lceil N\alpha\rceil+1)+i-1}^N\right]
\right\}.
\end{array}
$$
Using the fact that, given ${\cal G}_{q-1}^N$, the $\tilde X_q^{i}$'s, $i\in\{1,\dots,\lceil N\a\rceil\}$ are i.i.d. random vectors (see equation (\ref{decor-refbis}) page \pageref{decor-refbis}), and that with similar arguments they are independent of the subsample strictly below $L_q^N$, it is clear that $(Z^N_j)_{0\leq j\leq (p+1)(\lceil N\alpha\rceil+1)-1}$ is a triangular array of martingale increments adapted to the filtration ${\cal J}$. It is then straightforward to check that 
$$
{\cal M}_p^{N}={\cal M}_p^{N,1}+{\cal M}_p^{N,2}=\sum_{j=0}^{(p+1)(\lceil N\alpha\rceil+1)-1} Z^N_j,
$$
which is indeed a ${\cal J}$-martingale.\medskip 

Multiplying this large martingale by $\sqrt{N}$, we can use the CLT
theorem for martingales
page 171 of \cite{p}. The Lindeberg condition is obviously satisfied
since $f$ is assumed bounded, and the limits of the conditional
variances are specified by Propositions \ref{oisx} and \ref{oisw}. This terminates the proof of Theorem \ref{theorem1}.
\hfill $\blacksquare$

\paragraph{Remark}
This martingale decomposition may be found far from intuitive, but it highlights the contributions to the global error of both the empirical quantile, and the sample error. Moreover, it allows us to have a conditionally i.i.d. sample, and to use well known statistical properties of empirical quantiles.

\subsection{Proof of Corollary \ref{cor1}}
Concerning the proof of $(i)$, we just notice that
$$\sqrt{N}\left(\hat{E}-E\right)=\sqrt{N}\left(\hat{E}-E\right)\un_{\hat{n}=n}+\sqrt{N}\left(\hat{E}-E\right)\un_{\hat{n}\neq n}.$$
Then, for any $\varepsilon>0$, we have
$$\P\left(\left|\sqrt{N}\left(\hat{E}-E\right)\un_{\hat{n}\neq n}\right|>\varepsilon\right)\leq\P(\hat{n}\neq n).$$
Now, recall that, by Theorem \ref{kjsc}, $L_{n-1}^N$ and $L_{n}^N$ converge almost surely to $L_{n-1}$ and $L_{n}$, which ensures that $\hat{n}$ converges almost surely to $n$. As a consequence,
$$\sqrt{N}\left(\hat{E}-E\right)\un_{\hat{n}\neq n}\xrightarrow[N\to\infty]{\P}0.$$
Next, we have
$$\sqrt{N}\left(\hat{E}-E\right)\un_{\hat{n}=n}=\a^n\un_{\hat{n}=n}\times \sqrt{N}\left(\eta_n^N(f)-\eta_n(f)\right).$$
The first term on the right hand side converges in probability to $\a^n$ and, according to Theorem \ref{theorem1}, the second one converges in distribution to a Gaussian variable with variance $\Gamma(f)$. Putting all pieces together, we have shown that 
$$\sqrt{N}\left(\hat{E}-E\right)\xrightarrow[N\to\infty]{{\cal D}}{\cal N}(0,\a^{2n}\Gamma(f)).$$
Obviously, $(ii)$ is a direct application of this result with $f=\un_{S(\cdot)\geq L^\star}$. For $(iii)$, we have
\begin{align}
\sqrt{N}\left(\hat{C}-C\right)&=\sqrt{N}\left(\frac{\eta_n^N(f)}{\eta_n^N(\un_{S(\cdot)\geq L^\star})}-\frac{\eta_n(f)}{\eta_n(\un_{S(\cdot)\geq L^\star})}\right)\nonumber\\
&=\frac{\eta_n(\un_{S(\cdot)\geq L^\star})}{\eta_n^N(\un_{S(\cdot)\geq L^\star})}\times\sqrt{N}\left(\eta_n^N(g)-\eta_n(g)\right),\nonumber
\end{align}
where 
$$g=\frac{\un_{S(\cdot)\geq L^\star}}
{\eta_n(\un_{S(\cdot)\geq L^\star})}\left(f-\frac{\eta_n(f)}{\eta_n(\un_{S(\cdot)\geq L^\star})}\right).$$
Since $f=f\times\un_{S(\cdot)\geq L^\star}$, it is clear that $\e_n(g)=0$. Taking into account that $\eta_n(\un_{S(\cdot)\geq L^\star})=r$, we get
$$\Gamma(g)=\sum_{p=0}^n\e_p(\overline{Q}_{p,n}(g)^2)\quad\mbox{with}\quad g=\frac{\un_{S(\cdot)\geq L^\star}}{r}\left(f-\frac{\e_{n}(f)}{r}\right).$$ 
Moreover, we know from Theorem \ref{kjsc} that
$$\frac{\eta_n^N(\un_{S(\cdot)\geq L^\star})}{\eta_n(\un_{S(\cdot)\geq L^\star})}\xrightarrow[N\to\infty]{\P}1.$$
This concludes the proof of Corollary \ref{cor1}.

\hfill  $\blacksquare$

\section{Technical results}\label{technical}

This section gathers some general results which are used for establishing the proofs of Theorem \ref{kjsc} and Theorem \ref{theorem1}.

\subsection{Some regularity results}\label{sasc}
For $\mu$ an empirical or absolutely continuous probability distribution (like in Section \ref{prel}), and $K$ a transition kernel, we define the transition kernel $M_\mu$ as the truncated version of $K$ with respect to $\mu$, that is
$$M_\mu(x,dy) = G_\mu(x)(K(x,dy) G_\mu(y) + K(1-G_\mu)(x)\delta_x(dy)) + (1-G_\mu(x))\delta_x(dy).$$
Our first result is quite general but will be of constant use in the other proofs.

\begin{pro}\label{procont1}
Assume that $\nu(S^{-1}(\{L_\nu\}))=0$ and that $|L_\mu-L_\nu|\leq\delta$, then there exist two transition kernels $M^{\delta,-}$ and $M^{\delta,+}$ such that
\begin{itemize}
\item[(i)] $M^{\delta,-} \leq M_\mu \leq M^{\delta,+}$,
\item[(ii)] for all $f\in L^1(\nu)\cap L^1(\nu K)$, $\lim_{\delta \to 0} | \nu(M^{\delta,+} - M^{\delta,-})(f)| =0$.
\end{itemize}
Moreover, the same result holds if we replace respectively $M_\mu$ with $R_\mu=G_\mu M_\mu$, as well as $M^{\delta,-}$ with $R^{\delta,-}$, and  $M^{\delta,+}$ with $R^{\delta,+}$.
\end{pro}

Before proving this result, let us say something about the way we are going to apply it. Typically, we will consider the case where $\nu=\eta_p$ and $K=K_{p+1}$. Since $\eta_p\leq\alpha^{-p}\eta$ and recalling that $K_{p+1}$ is $\eta$ invariant, it is clear that if $f$ belongs to $L^1(\eta)$, then $f$ is in $L^1(\eta_p)\cap L^1(\eta_p K_{p+1})$ as well. Moreover, the absolute continuity of $\eta$ ensures that $\eta_p(S^{-1}(\{L_{\eta_p}\}))=0$.
\medskip 

{\bf Proof}\quad We will first prove the result for $M_\mu$, the other case is similar, just a bit simpler. We can decompose $M_\mu = M^0 + M^1 + M^2$ with
$$\left\{\begin{array}{l}
M^0(x,dy) = G_\mu(x) G_\mu(y) K(x,dy)\\
M^1(x,dy) = G_\mu(x) K(1-G_\mu)(x) \delta_x(dy)\\
M^2(x,dy) = (1-G_\mu(x))\delta_x(dy).\end{array}\right.$$
By construction, $G_{L_\nu+\delta}\leq G_\mu \leq G_{L_\nu-\delta}$. So we can take 
$$M^{\delta,+} = M^{0,\delta,+} +M^{1,\delta,+} +M^{2,\delta,+}$$ 
with
$$\left\{\begin{array}{l}
M^{0,\delta,+}(x,dy) = G_{L_\mu-\delta}(x) G_{L_\mu-\delta}(y) K(x,dy)\\
M^{1,\delta,+}(x,dy)= G_{L_\mu-\delta}(x) K(1-G_{L_\mu+\delta})(x) \delta_x(dy)\\
M^{2,\delta,+}(x,dy) = (1-G_{L_\mu+\delta}(x))\delta_x(dy)
\end{array}\right.$$
and similarly, we can take $M^{\delta,-} = M^{0,\delta,-} +M^{1,\delta,-} +M^{2,\delta,-}$ with
$$\left\{\begin{array}{l}
M^{0,\delta,-}(x,dy) = G_{L_\mu+\delta}(x) G_{L_\mu(y)+\delta}(y) K(x,dy)\\
M^{1,\delta,-}(x,dy)= G_{L_\mu+\delta}(x) K(1-G_{L_\mu-\delta})(x) \delta_x(dy)\\
M^{2,\delta,-}(x,dy) = (1-G_{L_\mu-\delta}(x))\delta_x(dy).
\end{array}\right.$$
Then $(i)$ is obviously satisfied. For $(ii)$, we clearly have for all $x\notin S^{-1}(\{L_\nu\})$, 
$$(M^{\delta,+}-M^{\delta,-})(f)(x)\xrightarrow[\delta\to 0]{} 0.$$
Moreover, a straightforward computation reveals that
$$\left|(M^{\delta,+}-M^{\delta,-})(f)\right|\leq K(|f|)+ 2|f|,$$ 
which belongs to $L^1(\nu)$ by assumption on $f$. We conclude using Lebesgue's dominated convergence theorem. For the other case, we can apply the same reasoning, by noticing that $R_\mu = G_\mu M_\mu = M^0 + M^1$, so that one can take 
$$R^{\delta,+} = M^{0,\delta,+} +M^{1,\delta,+}\hspace{1cm}\mbox{and}\hspace{1cm} R^{\delta,-} = M^{0,\delta,-} +M^{1,\delta,-}.$$
\hfill $\blacksquare$
\medskip

In the upcoming result, $(\nu_N)$ is a sequence of empirical probability measures on $\R^d$. We do not need to make further assumptions on its points for now.
Moreover, let $\nu$ be a fixed and absolutely continuous probability measure on $\R^d$. Denote respectively by $L$ and $L_N$ the $(1-\alpha)$ quantiles of $\nu$ and $\nu_N$ with respect to the mapping $S$ as defined in Section \ref{ad}, by ${\cal A}=\{x\in\R^d: S(x)\geq L\}$ and ${\cal A}_N=\{x\in\R^d: S(x)\geq L_N\}$ the associated level sets, and by $G(x)=\un_{{\cal A}}(x)$ and $G_N(x)=\un_{{\cal A}_N}(x)$ the related potential functions.
We will also assume that the probability measure $\nu \circ S^{-1}$ has a density, and that this density  is continuous and strictly 
positive at $L$.
\medskip

Moreover, if $K$ is a transition kernel on $\R^d$, we denote respectively by $M$ and $M_N$ its truncated versions according to $L$ and $L_N$, meaning that
$$M(x,dx')=\un_{\bar{{\cal A}}}(x)\d_x(dx')+\un_{{\cal A}}(x)(K(x,\bar{{\cal A}})\d_x(dx')+K(x,dx')\un_{{\cal A}}(x')),$$
and $M_N$ accordingly. The action of the mapping $\Pi$ on $\nu$ and $\nu_N$ is then defined as $\Pi(\nu)=\nu GM$ and $\Pi(\nu_N)=\nu_NG_NM_N$. The following result exhibits the continuity of $\Pi$.

\begin{pro}\label{asc}
With the previous notation , if for any $f\in L^1(\nu)\cap L^1(\nu K)$, one has
$$\nu_N(f)\xrightarrow[N\to\infty]{}\nu(f) \mbox{ a.s.\!\! (resp.\! in probability) }$$
then
\begin{enumerate}
\item[](i) $L_N\xrightarrow[N\to\infty]{}L$ a.s.\!\! (resp.\! in probability).
\item[](ii) $\Pi(\nu_N)(f)\xrightarrow[N\to\infty]{}\Pi(\nu)(f)$ a.s.\!\! (resp.\! in probability).
\end{enumerate}
\end{pro}

{\bf Proof}\quad 
We prove only the convergence a.s., the convergence in probability will follow using a.s. convergence of subsequences. 
\medskip

To prove $(i)$, let us fix $\varepsilon>0$ and let us denote by $F$ the cdf of the absolutely continuous probability measure $\nu\circ S^{-1}$. By assumption on $F$, there exist two strictly positive real numbers $\d^-$ and $\d^+$ such that
$$F(L-\varepsilon)=1-\a-\d^-\quad\quad\mbox{and}\quad\quad F(L+\varepsilon)=1-\a+\d^+.$$   
Applying the almost sure convergence of $\nu_N(f)$ to $\nu(f)$ respectively with $f=\un_{S(\cdot)\leq L-\varepsilon}$ and $f=\un_{S(\cdot)\leq L+\varepsilon}$, we get that for $N$ large enough,
$$\nu_N(\un_{S(\cdot)\leq L-\varepsilon})\leq 1-\a-\frac{\d^-}{2}\quad\quad\mbox{and}\quad\quad \nu_N(\un_{S(\cdot)\leq L+\varepsilon})\geq 1-\a+\frac{\d^+}{2}.$$
This ensures that, for $N$ large enough, $|L_N-L|\leq\varepsilon$. Since $\varepsilon$ is arbitrary, point $(i)$ is proved.
\medskip

Now we prove $(ii)$. From $(i)$, for any $\delta>0$, for $N$ larger than some random $N_0$, we have that $|L_{\nu_N}-L_\nu|\leq\delta$ and we are in a position to apply Proposition \ref{procont1}. Moreover, the triangular inequality gives
\begin{eqnarray*}
|(\Pi(\nu_N)-\Pi(\nu))(f)| &=& |(\nu_N G_{N} M_{N} - \nu G M) (f)|\\
&\leq& |\nu_N ( G_{N} M_{N} - G M) (f)|+  |(\nu_N-\nu)(G M(f))|
\end{eqnarray*}
where the second term can be made arbitrarily small by assumption. 
For the first term, we have
$$
 |\nu_N ( G M - G M) (f)|
 \leq  \nu_N (  R^{\delta,+} - R^{\delta,-}) (|f|),
 $$
 which converges to $|\nu (  R^{\delta,+} - R^{\delta,-}) (f)|$ by assumption. We conclude by choosing $\delta$ such that the limit is arbitrarily small.
 \hfill  $\blacksquare$
\medskip

Our next result will be used in the proof of Proposition \ref{VieillePuteAlbanaise}.

\begin{cor}\label{tilde}
For any $q\in\{0,\dots,n-1\}$, for all $f\in L^2(\eta)$,
$$\tilde\eta_q^N(f)\xrightarrow[N\to\infty]{\P}\eta_{q+1}(f),$$
and for all $f\in \Ba(\R^d)$,
$$\tilde\eta_q^N(f)\xrightarrow[N\to\infty]{a.s.}\eta_{q+1}(f).$$
\end{cor}

{\bf Proof}\quad We only treat the case where $f$ belongs to $L^2(\eta)$. By (\ref{decor-ref}), we have
$$\tilde\eta_q^N(f)=\Psi_{G_{\eta^N_q}}(\eta^N_q)(f)=\frac{N}{\lceil N\alpha\rceil}\eta_q^N(G_{\eta^N_q}\times f).$$
Assume that the transition kernel $K_{q+1}$ is the identity, that is $K_{q+1}(x,\cdot)=\delta_x$, then by (\ref{io}) and the definition of $\Pi_{q+1}$, we may write
$$\tilde\eta_q^N(f)=\frac{N}{\lceil N\alpha\rceil}\Pi_{q+1}(\eta_q^N)(f).$$
From Theorem \ref{kjsc}, we know that
$$\eta_q^N(f)\xrightarrow[N\to\infty]{\P}\eta_q(f).$$
Thus, since 
$$L^2(\eta)\subset L^1(\eta_q)=L^1(\eta_q)\cap L^1(\eta_q K_{q+1}),$$ 
Proposition \ref{asc} yields
$$\tilde\eta_q^N(f)=\frac{N}{\lceil N\alpha\rceil}\Pi_{q+1}(\eta_q^N)(f)\xrightarrow[N\to\infty]{\P}\frac{1}{\alpha}\Pi_{q+1}(\eta_q)(f)=\frac{1}{\alpha}\eta_q(G_{q}\times f)=\eta_{q+1}(f).$$ 
\hfill  $\blacksquare$
\medskip

The upcoming corollary is at the core of the proofs of Propositions \ref{oisx} and \ref{oisw}.

\begin{cor}\label{aozuh}
For any $1\leq q\leq p<n$, any $f\in\Ba(\R^d)$, we have
$$\Psi_{G_{\eta^N_q}}\left(\Phi_q(\eta^N_{q-1})\right)(f)\xrightarrow[N\to\infty]{a.s.}\eta_{q+1}(f),$$
and for any $\beta>0$,
$$\Psi_{G_{\eta^N_q}}\left(\Phi_q(\eta^N_{q-1})\right)\left\{
\left(\left[\widetilde{Q}_{q,p,\eta^N_q}-\widetilde{Q}_{q,p}\right](f)\right)^\beta\right\}\xrightarrow[N\to\infty]{a.s.}0,$$
\end{cor}

{\bf Proof}\quad By Theorem \ref{kjsc}, we know that for all $f\in\Ba(\R^d)$, we have
$$\eta_{q-1}^N(f)\xrightarrow[N\to\infty]{a.s.}\eta_{q-1}(f).$$
Hence, by Proposition \ref{asc}, we deduce that for all $f\in\Ba(\R^d)$,
$$\Pi_q(\eta_{q-1}^N)(f)\xrightarrow[N\to\infty]{a.s.}\Pi_q(\eta_{q-1})(f).$$
Next, by (\ref{decor-reftetra}), we may write
$$\Psi_{G_{\eta^N_q}}\left(\Phi_q(\eta^N_{q-1})\right)(f)=\frac{\Pi_q(\eta^N_{q-1})(G_{\eta^N_q}f)}{\Pi_q(\eta^N_{q-1})(G_{\eta^N_q})}.$$
Still by Theorem \ref{kjsc}, we know that
$$L_{\eta_q^N}=L_q^N\xrightarrow[N\to\infty]{a.s.}L_q.$$
Thus, for any $\delta>0$, almost surely for $N$ large enough, one has 
$$G_{L_q+\delta}\leq G_{\eta^N_q},G_q\leq G_{L_q-\delta}$$
and the same reasoning as in the proof of Proposition \ref{asc} shows that for all $f\in\Ba(\R^d)$,  
$$\Psi_{G_{\eta^N_q}}\left(\Phi_q(\eta^N_{q-1})\right)(f)=\frac{\Pi_q(\eta^N_{q-1})(G_{\eta^N_q}f)}{\Pi_q(\eta^N_{q-1})(G_{\eta^N_q}).}\xrightarrow[N\to\infty]{a.s.}\frac{\Pi_q(\eta_{q-1})(G_{q}f)}{\Pi_q(\eta_{q-1})(G_{q})}=\eta_{q+1}(f).$$
For the second point, first notice that
$$\left[\widetilde{Q}_{q,p,\eta^N_q}-\widetilde{Q}_{q,p}\right](f)=[M_{q+1,\eta_q^N}-M_{q+1}](Q_{q+1,p}(f)).$$
Then, by the first point of Proposition \ref{procont1}, we deduce that almost surely for $N$ large enough,
$$\left|[M_{q+1,\eta_q^N}-M_{q+1}](Q_{q+1,p}(f))\right|\leq \left|[M_{q+1}^{\delta,+}-M_{q+1}^{\delta,-}](Q_{q+1,p}(f))\right|.$$
Therefore, by the previous point,
\begin{align}
&\limsup_{N\to\infty}\left|\Psi_{G_{\eta^N_q}}\left(\Phi_q(\eta^N_{q-1})\right)\left\{
\left(\left[\widetilde{Q}_{q,p,\eta^N_q}-\widetilde{Q}_{q,p}\right](f)\right)^\beta\right\}\right|\nonumber\\
&\leq\eta_{q+1}\left(\left|[M_{q+1}^{\delta,+}-M_{q+1}^{\delta,-}](Q_{q+1,p}(f))\right|^\beta\right).\nonumber
\end{align}
Finally, the desired result is just a consequence of the second point of Proposition \ref{procont1}. 
 \hfill  $\blacksquare$
\medskip

Basically, the previous results focused on the continuity of the operator $\Pi$. In the remainder of this subsection, we go one step further as we are interested in asymptotic expansions. We recall that 
$${\cal B}^a_q=\left\{g:\R^d\to\R,\ \exists (g_0\dots,g_{q-1})\in\Ba(\R^d)^q,\ g=\ck^a_1(g_0)\cdots \ck^a_q(g_{q-1})\right\},$$
and for $g\in{\cal B}^a_q$, $x\in\R^d$ and $L\in\R$, we denote 
$$H_q^{g,a}(x,L)=\int_{S(x')=L}\ck^a_{q+1}(x,x')g(x')\frac{\dbar x'}{|DS(x')|}.$$
Let us first generalize the notations of Assumption [${\cal H}^a$] to any probability measure $\nu$. As before, we typically have in mind the case where $\nu=\eta_q$ is the restriction of $\eta$ above level $L_{q-1}$, in which case Assumption $[{\cal H}_\nu^a]$ will be equivalent to Assumption [${\cal H}^a$]. If we consider the kernel $K_q$ instead of $\ck^a_q$, we will have exactly the same results, as it is a special case for which $a(x,x^\prime)$ is constant equal to $1$.
\medskip

\textbf{Assumption [${\cal H}_\nu^a$]}
\begin{enumerate}
\item[$(i)$] For any $q\geq 0$, the mapping $x\mapsto H_q^{1,a}(x,L_q)$ belongs to $L^2(\nu)$, that is 
$$\int\nu(dx)\left(\int_{S(x')=L_q}\ck^a_{q+1}(x,x')\frac{\dbar x'}{|DS(x')|}\right)^2<\infty,$$
and $\nu(H_q^{1,a}(.,L_q))>0$.
\item[$(ii)$] For any $q>0$, for any $g\in{\cal B}^a_q$, there exists $h\in L^2(\nu)$ such that for any $\varepsilon>0$, there exists $\delta>0$ such that for any $L\in[L_q-\delta,L_q+\delta]$ and for almost every $x\in\R^d$,
$$\left|H_q^{g,a}(x,L)-H_q^{g,a}(x,L_q)\right|\leq \varepsilon h(x).$$
\end{enumerate}
\medskip

The following result will be of constant use in the proof of Proposition \ref{VieillePuteAlbanaise}.

\begin{lem}\label{undeplus}
Assume that for any $f\in L^2(\nu)$, one has
$$\nu_N(f)\xrightarrow[N\to\infty]{\P}\nu(f).$$
Then, for any $g\in{\cal B}^a_q$ and any $\varphi\in\Ba(\R^d)$, under Assumption $[{\cal H}_\nu^a]$, one has
$$\nu_N\left(\varphi\int_{L_q}^{L_q^N}H_q^{g,a}(\cdot,L)dL\right)=(L_q^N-L_q)\nu(\varphi H_q^{g,a}(\cdot,L_q))+o_p(L_q^N-L_q).$$
\end{lem}

{\bf Proof}\quad We first choose $\varepsilon>0$. By point $(ii)$ of Assumption $[{\cal H}_\nu^a]$, the mapping $L\mapsto H_q^g(x,L)$ is continuous in the neighborhood of $L_q$ for $\nu$ almost every $x$. We consider $N$ large enough such that $L_q^N\in (L_q-\delta,L_q+\delta)$ with arbitrarily large probability, say $1-\gamma$. Hence, by the mean value theorem, there exists $\tilde{L}$ between $L_q$ and $L_q^N$ such that
$$\int_{L_q}^{L_q^N}H_q^{g,a}(x,L)dL=(L_q^N-L_q)\times H_q^{g,a}(x,\tilde{L}).$$ 
As a consequence,
$$\frac{\nu_N\left(\varphi\displaystyle\int_{L_q}^{L_q^N}H_q^{g,a}(\cdot,L)dL\right)}{L_q^N-L_q}=\nu_N(\varphi H_q^{g,a}(\cdot,L_q))+\nu_N(\varphi (H_q^{g,a}(\cdot,\tilde{L})-H_q^{g,a}(\cdot,L_q))).$$
Since $\varphi$ and $g$ are both bounded, point $(i)$ of Assumption $[{\cal H}_\nu^a]$ ensures that the function $\varphi H_q^{g,a}(\cdot,L_q)$ is in $L^2(\nu)$, so that by the hypothesis of Lemma \ref{undeplus}, 
$$\nu_N(\varphi H_q^{g,a}(\cdot,L_q))\xrightarrow[N\to\infty]{\P}\nu(\varphi H_q^{g,a}(\cdot,L_q)).$$
Furthermore, by point $(ii)$ of Assumption $[{\cal H}_\nu^a]$, we have
$$\left|\nu_N(\varphi (H_q^{g,a}(\cdot,\tilde{L})-H_q^{g,a}(\cdot,L_q)))\right|\leq  (\|\varphi\|\times\nu_N(h))\times \varepsilon,$$
where, since $h$ belongs to $L^2(\nu)$,
$$\nu_N(h)\xrightarrow[N\to\infty]{\P}\nu(h)<\infty.$$
Since $\varepsilon$ and $\gamma$ are arbitrary, the proof is complete.
\hfill  $\blacksquare$ 

\subsection{Proof of Proposition \ref{oisx}}
Our goal is to prove that, for any $q\leq p$ and any $f\in\Ba(\R^d)$, we have
$$
N\a^{2(q-p)}~\V\left(\left.\Pi_{q,p}(\eta_q^N)(f)\right|{\cal G}_{q-1}^N\right)\xrightarrow[N\to\infty]{a.s.}\e_q(\overline{Q}_{q,p}(f)^2)-\a^{-1}~\e_p(f)^2.\nonumber
$$ 
By (\ref{bla}), we have
\begin{eqnarray*}
\Pi_{q,p}(\eta^N_q)&=&\frac{\lceil N\alpha\rceil}{N}~~\widetilde{\eta}^N_q\widetilde{Q}_{q,p,\eta^N_q}
\end{eqnarray*}
with the measure $\widetilde{\eta}^N_q$ defined in (\ref{decor-ref}).
This shows that
\begin{align}
N\a^{2(q-p)}\V\left(\left.\Pi_{q,p}(\eta_q^N)(f)\right|{\cal G}_{q-1}^N\right)=&\left(\frac{\lceil N\alpha\rceil}{N}\right)^2N\a^{2(q-p)}\V\left(\left.\widetilde{\eta}^N_q\widetilde{Q}_{q,p,\eta^N_q}(f)\right|{\cal G}_{q-1}^N\right)\nonumber\\
=&\ \a^{2(q-p)}~\frac{\lceil N\alpha\rceil}{N}~\V\left(\left.\widetilde{Q}_{q,p,\eta^N_q}(f)(\widetilde{X}_q^{1})\right|{\cal G}_{q-1}^N\right).\label{ref-VV}
\end{align}
On the other hand, we have, thanks to (\ref{decor-refbis}), 
\begin{align}
&\V\left(\left.\left[\widetilde{Q}_{q,p,\eta^N_q}-\widetilde{Q}_{q,p}\right](f)(\widetilde{X}^{1}_q)\right|{\cal G}_{q-1}^N\right)\nonumber\\
&=\Psi_{G_{\eta^N_q}}\left(\Phi_q(\eta^N_{q-1})\right)\left\{
\left(\left[\widetilde{Q}_{q,p,\eta^N_q}-\widetilde{Q}_{q,p}\right](f)\right)^2\right\}\nonumber\\
&\ \ \ -\left( \Psi_{G_{\eta^N_q}}\left(\Phi_q(\eta^N_{q-1})\right)\left(\left[\widetilde{Q}_{q,p,\eta^N_q}-\widetilde{Q}_{q,p}\right](f)\right)\right)^2.\nonumber
\end{align}
By the second point of Corollary \ref{aozuh}, we deduce that
$$\V\left(\left.\left[\widetilde{Q}_{q,p,\eta^N_q}-\widetilde{Q}_{q,p}\right](f)(\widetilde{X}^{1}_q)\right|{\cal G}_{q-1}^N\right)\xrightarrow[N\to\infty]{a.s.}0.$$
In other words, coming back to (\ref{ref-VV}) and applying the first point of Corollary \ref{aozuh}, we have obtained
\begin{align}
&N\a^{2(q-p)}~\V\left(\left.\Pi_{q,p}(\eta_q^N)(f)\right|{\cal G}_{q-1}^N\right)\nonumber\\
&\xrightarrow[N\to\infty]{a.s.}\a^{2(q-p)+1}\left\{\e_{q+1}([M_{q+1}Q_{q+1,p}(f)]^2)-(\e_{q+1}M_{q+1}Q_{q+1,p}(f))^2\right\}.\nonumber
\end{align}
Using elementary computations, it is easy to check that
\begin{align}
&\a^{2(q-p)+1}~\left\{\e_{q+1}([M_{q+1}Q_{q+1,p}(f)]^2)-(\e_{q+1}Q_{q+1,p}(f))^2\right\}\nonumber\\
&=\e_{q}(\overline{Q}_{q,p}(f)^2)-\a^{-1}~\e_p(f)^2,\nonumber
\end{align}
which terminates the proof of Proposition \ref{oisx}.
\hfill $\blacksquare$

\subsection{Proof of Proposition \ref{oisw}}
We intend to show that, for any $q\leq p$ and any $f\in\Ba(\R^d)$, we have
$$
\E\left[\left.\sqrt{N}\alpha^{q-p}~\epsilon_{q}^N~\E\left[\left.\Pi_{q,p}(\eta_q^N)(f)\right|{\cal G}_{q-1}^N\right]\right| {\cal F}_{q-1}^N\right]\xrightarrow[N\to\infty]{\P} 0,
$$
and 
$$
\V\left[\left.\sqrt{N}\alpha^{q-p}~\epsilon_{q}^N~\E\left[\left.\Pi_{q,p}(\eta_q^N)(f)\right|{\cal G}_{q-1}^N\right]\right| {\cal F}_{q-1}^N\right]
 \xrightarrow[N\to\infty]{\P} \frac{1-\alpha}{\alpha}\ \eta_p(f)^2.
$$
The proof is carried out given ${\cal F}_{q-1}^N$. We begin like in the proof of Proposition \ref{oisx}. From (\ref{ref-cond-G}) and the definition of $\rho_N$, recall that
$$\E\left[\Pi_{q,p}(\eta^N_q)(f)|{\cal G}_{q-1}^N\right]=\alpha\rho_N~\Psi_{G_{\eta^N_q}}\left(\Phi_q(\eta^N_{q-1})\right)\left(\widetilde{Q}_{q,p,\eta^N_q}(f)\right).
$$ 
Hence, the quantity of interest in Proposition \ref{oisw} may be rewritten as follows
\begin{align}
&\sqrt{N}\alpha^{q-p}~\epsilon^N_q~\E\left[\left.\Pi_{q,p}(\eta_q^N)(f)\right|{\cal G}_{q-1}^N\right]\nonumber\\
&=\left[\sqrt{N}~\alpha~\epsilon^N_q\right]~\rho_N~\alpha^{q-p}~\Psi_{G_{\eta^N_q}}\left(\Phi_q(\eta^N_{q-1})\right)\left(\widetilde{Q}_{q,p,\eta^N_q}(f)\right). \nonumber 
\end{align}
We are going to prove first that 
\begin{equation*}
\rho_N\Psi_{G_{\eta^N_q}}\left(\Phi_q(\eta^N_{q-1})\right)\left(\widetilde{Q}_{q,p,\eta^N_q}(f)\right)\xrightarrow[N\to\infty]{a.s.}\alpha^{p-q-1}~\eta_p(f).
\end{equation*}
Because  $|\rho_N-1|\leq 1/(N\alpha)$, the factor $\rho_N$ is unimportant.
 As in the proof of Proposition \ref{oisx}, we consider the decomposition
$$
\begin{array}{l}
\Psi_{G_{\eta^N_q}}\left(\Phi_q(\eta^N_{q-1})\right)\left(\widetilde{Q}_{q,p,\eta^N_q}(f)\right)\\
\\
=\Psi_{G_{\eta^N_q}}\left(\Phi_q(\eta^N_{q-1})\right)\left(\widetilde{Q}_{q,p,\eta^N_q}(f)-\widetilde{Q}_{q,p}(f)\right)+\Psi_{G_{\eta^N_q}}\left(\Phi_q(\eta^N_{q-1})\right)\left(\widetilde{Q}_{q,p}(f)\right).
\end{array}$$
The second point of Corollary \ref{aozuh} implies that
$$\left|\Psi_{G_{\eta^N_q}}\left(\Phi_q(\eta^N_{q-1})\right)\left(\left[\widetilde{Q}_{q,p,\eta^N_q}-\widetilde{Q}_{q,p}\right](f)\right)\right|\xrightarrow[N\to\infty]{a.s.}0,$$
while the first point of Corollary \ref{aozuh} ensures that
$$\Psi_{G_{\eta^N_q}}\left(\Phi_q(\eta^N_{q-1})\right)\left(\tilde{Q}_{q,p}(f)\right)\xrightarrow[N\to\infty]{a.s.}\e_{q+1}\left(\tilde{Q}_{q,p}(f)\right)=\alpha^{p-q-1}\eta_p(f).$$
As $f$ is bounded, by Lebesgue's dominated convergence theorem, the above convergence also holds in $L^2$. 

\medskip

Now we prove the first assertion of Proposition \ref{oisw}. Let us denote for a moment 
$$Z_N=\rho_N \Psi_{G_{\eta^N_q}}\left(\Phi_q(\eta^N_{q-1})\right)\left(\widetilde{Q}_{q,p,\eta^N_q}(f)\right)$$ 
and by $z=\alpha^{p-q-1}~\eta_p(f)$ its deterministic limit. We have just shown that $Z_N-z$ converges to $0$ in $L^2$. This implies that  
$$\E\left[(Z_N-z)^2\ |\ \F_{q-1}^N\right]\xrightarrow[N\to\infty]{\P}0.$$
By very similar arguments, we can also see that 
$$\E\left[(Z_N^2-z^2)^2\ |\ \F_{q-1}^N\right]\xrightarrow[N\to\infty]{\P}0.$$  
We have
$$
\sqrt{N}\E\left[\epsilon_q^N Z_N\ |\ \F_{q-1}^N\right]=\sqrt{N}\E\left[\epsilon_q^N(Z_N-z)\ |\ \F_{q-1}^N\right] + \sqrt{N}z\E\left[\epsilon_q^N\ |\ \F_{q-1}^N\right].
$$
Since $|\rho_N-1|\leq 1/(N\alpha)$, we have
$$\epsilon_q^N=\frac{1}{\alpha}\left(\alpha-\Phi_q(\eta_{q-1}^N)(G_{\eta^N_{q}})\right)+o(1/\sqrt{N}),$$
and the convergence of the second term is a direct consequence of the first result of Lemma~\ref{lem.TCL.1} below. For the first term, we apply Cauchy-Schwarz inequality
$$
\left| \sqrt{N}\E\left[\epsilon_q^N(Z_N-z)\ |\ \F_{q-1}^N\right] \right|
 \leq
 \sqrt{N\E\left[(\epsilon_q^N)^2\ |\ \F_{q-1}^N\right]} \sqrt{\E\left[(Z_N-z)^2\ |\ \F_{q-1}^N\right]},
 $$
 which converges in probability to $0$ by the second result of Lemma~\ref{lem.TCL.1}, and the $L^2$ convergence of $Z_N-z$.
 
\medskip
For the second assertion of Proposition \ref{oisw}, we write
$$
{N}\E\left[(\epsilon_q^N Z_N)^2\ |\ \F_{q-1}^N\right]={N}\E\left[(\epsilon_q^N)^2(Z_N^2-z^2)\ |\ \F_{q-1}^N\right] + {N}z^2\E\left[(\epsilon_q^N)^2\ |\ \F_{q-1}^N\right].
$$
The convergence of the second term is a direct consequence of the second result of Lemma~\ref{lem.TCL.1}. For the first term we use Cauchy--Schwarz again
$$
 {N}\E\left[(\epsilon_q^N)^2(Z_N^2-z^2)\ |\ \F_{q-1}^N\right] 
 \leq
 \sqrt{N^2\E\left[(\epsilon_q^N)^4\ |\ \F_{q-1}^N\right]} \sqrt{\E\left[(Z_N^2-z^2)^2\ |\ \F_{q-1}^N\right]}
 $$
 and we conclude similarly, using the third result of Lemma~\ref{lem.TCL.1}.
 \hfill  $\blacksquare$ 

\begin{lem}
\label{lem.TCL.1}
For any integer $q$, we have
$$
\sqrt{N}~\E\left[\left.\Phi_q(\eta_{q-1}^N)(G_{\eta^N_{q}})-\alpha\ \right| \ \F_{q-1}^N \right]\xrightarrow[N\to\infty]{\P} 0,
$$
$$
N~\E\left[\left.\left(\Phi_q(\eta_{q-1}^N)(G_{\eta^N_{q}})-\alpha\right)^2\ \right| \ \F_{q-1}^N \right] \xrightarrow[N\to\infty]{\P} \alpha(1-\alpha),
$$
and
$$
N^2~\E\left[\left.\left(\Phi_q(\eta_{q-1}^N)(G_{\eta^N_{q}})-\alpha\right)^4\ \right| \ \F_{q-1}^N \right]={\cal O}_p(1).
$$
\end{lem} 
 
\paragraph{Proof}

Here again, the reasoning is made given ${\cal F}^N_{q-1}$. Recall that $(X_q^i)_{1\leq i\leq N}$ is an i.i.d. sample with common law $\Phi_{q }\left(\eta_{q-1}^N\right)$. Accordingly, let us denote $(Y_q^i)_{1\leq i\leq N}=(S(X_q^i))_{1\leq i\leq N}$. Also, for any real number $L$, define the function 
$$F_N(L)=1-\Phi_q(\eta_{q-1}^N)(G_L),$$ 
which is more or less a cumulative distribution function. The function $F_N$ is continuous except at a finite number of values, namely at most the $\lceil\alpha N\rceil$ largest values among the $Y_{q-1}^i$'s.
\medskip

Starting from the sample $(Y_q^i)_{1\leq i\leq N}$, we also construct a new sample ${\bf U}=(U_q^i)_{1\leq i\leq N}$ as follows. If $Y_q^i$ is a point of continuity of $F_N$, then $U_q^i=F_N(Y_q^i)$, otherwise we draw $U_k^i$
uniformly in the interval 
$$(F_N(Y_q^i),\lim_{h\to 0^+}F_N(Y_q^i+h)).$$ 
It is then a simple exercise (see for example \cite{shorack86a}, page 102) to
check that ${\bf U}=(U_q^1,\dots,U_q^N)$ is an i.i.d. sample with distribution
${\cal U}(0,1)$. Denoting 
$$U^N_{1-\alpha}:=U^N_{(\frac{k_N}{N})}=U^N_{(\frac{\lceil(1-\alpha)N\rceil}{N})},$$ 
we may write
\begin{align}
\Phi_q(\eta_{q-1}^N)(G_{\eta^N_{q}})-\alpha &=(U^N_{1-\alpha}-(1-\Phi_q(\eta_{q-1}^N)(G_{\eta^N_{q}})))+((1-\alpha)-U^N_{1-\alpha})\nonumber\\
&=(U^N_{1-\alpha}-F_N(L_q^N))+((1-\alpha)-U^N_{1-\alpha}).\label{eqf1}
\end{align}
The first term can easily be bounded in absolute value thanks to the following lemma, whose proof is detailed in Section \ref{secmulti}.

\begin{lem}
\label{lemmulti}
For any integer $q$ and any $\ell\in\{1,2,4\}$, we have
$$
N^{\ell/2}~\E\left[\left.\left(U^N_{1-\alpha}-F_N(L_q^N)\right)^\ell\ \right| \ \F_{q-1}^N \right] \xrightarrow[N\to\infty]{\P} 0.
$$
\end{lem}

For the second term in  (\ref{eqf1}), we have 
\begin{align}
(1-\alpha)-U^N_{1-\alpha} &= \left((1-\alpha)-\frac{k_N}{N}\right)+\left(\frac{k_N}{N}-U^N_{(\frac{k_N}{N})}\right).\label{eqf2}
\end{align}
The first term is deterministic and goes to $0$. For the second term, it is well known that (see e.g. \cite{shorack86a} page 97)
$$U^N_{(\frac{k_N}{N})}\sim{\rm Beta}(k_N,N-k_N+1).$$ 
Therefore 
$$\E\left[\left.U^N_{(\frac{k_N}{N})}\right| \ \F_{q-1}^N\right]=\frac{k_N}{N}\hspace{1cm}\mbox{and}\hspace{1cm}\V\left(\left.U^N_{(\frac{k_N}{N})}\right| \ \F_{q-1}^N\right)=\frac{k_N(N-k_N+1)}{(N+1)^2(N+2)},$$ 
so that 
$$N\times\V\left(\left.U^N_{(\frac{k_N}{N})}\right| \ \F_{q-1}^N\right)=N\times\frac{k_N(N-k_N+1)}{(N+1)^2(N+2)}\xrightarrow[N\to\infty]{} \alpha(1-\alpha).$$
We also have
\begin{align}
&N^2\times\E\left[\left.\left(U^N_{(\frac{k_N}{N})}-\frac{k_N}{N}\right)^4\right| \ \F_{q-1}^N\right]\nonumber\\
&\quad\quad=N^2\times\frac{3 k_N (N-k_N+1) (2(N+1)^2 + k_N (N-k_N+1) (N-5))}{(N+1)^4(N+2)(N+3)(N+4)},\nonumber
\end{align}
which is obviously bounded.\medskip

Let us prove the first assertion of Lemma \ref{lem.TCL.1}. From (\ref{eqf1}) and (\ref{eqf2}) we deduce
\begin{align}
&\sqrt{N}~\E\left[\left.\Phi_q(\eta_{q-1}^N)(G_{\eta^N_{q}})-\alpha\ \right| \ \F_{q-1}^N \right]\nonumber\\
&=\sqrt{N}~\E\left[\left.U^N_{1-\alpha}-F_N(L_q^N) \right| \ \F_{q-1}^N \right]+\sqrt{N}\left((1-\alpha)-\frac{k_N}{N}\right).\nonumber
\end{align}
The first term goes to 0 in probability thanks to Lemma \ref{lemmulti}, the second one is deterministic and goes to 0 since $|(1-\alpha)-k_N/N|\leq 1/N$.\medskip   

For the second assertion of Lemma \ref{lem.TCL.1}, relation (\ref{eqf1}) gives
\begin{align}
&N~\E\left[\left.\left(\Phi_q(\eta_{q-1}^N)(G_{\eta^N_{q}})-\alpha\right)^2\ \right| \ \F_{q-1}^N \right]\nonumber\\
&=N~\E\left[\left.\left(U^N_{1-\alpha}-F_N(L_q^N)\right)^2 \right| \ \F_{q-1}^N \right]+N~\E\left[\left.\left((1-\alpha)-U^N_{1-\alpha}\right)^2 \right| \ \F_{q-1}^N \right]\nonumber\\
&\quad+N~\E\left[\left.\left(U^N_{1-\alpha}-F_N(L_q^N)\right)\left((1-\alpha)-U^N_{1-\alpha}\right) \right| \ \F_{q-1}^N \right].\label{eqf3}
\end{align}
Here again, the first term goes to 0 in probability thanks to Lemma \ref{lemmulti}. For the second one, just notice that
\begin{align}
&\E\left[\left.\left((1-\alpha)-U^N_{1-\alpha}\right)^2 \right| \ \F_{q-1}^N \right]=\V\left(\left.U^N_{(\frac{k_N}{N})}\right| \ \F_{q-1}^N\right)+\left((1-\alpha)-\frac{k_N}{N}\right)^2,\nonumber
\end{align}
which implies that
$$N~\E\left[\left.\left((1-\alpha)-U^N_{1-\alpha}\right)^2 \right| \ \F_{q-1}^N \right]\xrightarrow[N\to\infty]{}\alpha(1-\alpha).$$
Finally, Cauchy-Schwarz shows that the last term in (\ref{eqf3}) goes to 0 in probability, and the second assertion of Lemma \ref{lem.TCL.1} is established.\medskip 

Concerning the third assertion of Lemma \ref{lem.TCL.1}, it suffices to remark that
\begin{align}
&N^2~\E\left[\left.\left(\Phi_q(\eta_{q-1}^N)(G_{\eta^N_{q}})-\alpha\right)^4\ \right| \ \F_{q-1}^N \right]\nonumber\\
&\leq 4N^2~\E\left[\left.\left(U^N_{1-\alpha}-F_N(L_q^N)\right)^4 \right| \ \F_{q-1}^N \right]+4N^2~\E\left[\left.\left((1-\alpha)-U^N_{1-\alpha}\right)^4 \right| \ \F_{q-1}^N \right].\nonumber
\end{align}
Then Lemma \ref{lemmulti} and the fourth moment of the Beta distribution of interest ensure that this quantity is bounded in probability. This completes the proof of Lemma \ref{lem.TCL.1}. 
\hfill  $\blacksquare$

\subsection{Proof of Proposition \ref{VieillePuteAlbanaise}}
We have to show that, under Assumption $[{\cal H}]$, for any $q\leq p$ and any $f\in\Ba(\R^d)$ such that $f=f\times\un_{S(\cdot)\geq L^\star}$, 
$$\sqrt{N}~\left(\alpha~(1-\epsilon_{q}^N)\E\left[\left.\Pi_{q,p}(\eta_q^N)(f)\right|{\cal G}_{q-1}^N\right]-\Pi_{q-1,p}(\eta_{q-1}^N)(f)\right)\xrightarrow[N\to\infty]{\P}0.$$
It turns out that the proof is quite technical and requires several auxiliary results whose proofs are postponed to the end of the present section. Here again, the reasoning is carried out given ${\cal F}_{q-1}^N$. By (\ref{yoyo}), (\ref{io}) and the definition of $\rho_N$, we have
\begin{eqnarray*}
\Pi_{q-1,p}(\eta^N_{q-1})(f)&=&\eta^N_{q-1}(G_{\eta^N_{q-1}})\Phi_{q}(\eta^N_{q-1})Q_{q,p}(f)\\
&=&\alpha\rho_N~
\Phi_{q}(\eta^N_{q-1})Q_{q,p}(f)
\end{eqnarray*}
and by (\ref{yo}) and (\ref{ref-cond-G}),
$$
\begin{array}{l}
\alpha~(1-\epsilon_{q}^N)\E\left[\left.\Pi_{q,p}(\eta_q^N)(f)\right|{\cal G}_{q-1}^N\right]\\
\\
=\rho_N~\Phi_{q }(\eta_{q-1}^N)(G_{\eta^N_q})~
\E\left[\Pi_{q,p}(\eta^N_q)(f)~|~{\cal G}_{q-1}^N\right]\\
\\
=\alpha\rho_N^2~\Phi_{q }(\eta_{q-1}^N)(G_{\eta^N_q})~\Psi_{G_{\eta^N_q}}\left(\Phi_q(\eta^N_{q-1})\right)\widetilde{Q}_{q,p,\eta^N_q}(f)\\
\\
=\alpha\rho_N^2~\Phi_{q }(\eta_{q-1}^N)\left(Q_{q,p,\eta^N_q}(f)\right).
\end{array}$$
Since $f$ is bounded and $\rho_N-1={\cal O}(N^{-1})$, this implies that
$$
\begin{array}{l}
\alpha~(1-\epsilon_{q}^N)\E\left[\left.\Pi_{q,p}(\eta_q^N)(f)\right|{\cal G}_{q-1}^N\right]-\Pi_{q-1,p}(\eta_{q-1}^N)(f)\\
\\
=\alpha\rho_N~\Phi_{q }\left(\eta_{q-1}^N\right)\left(\left[Q_{q,p,\eta^N_q}-Q_{q,p}\right](f)\right)+{\cal O}(N^{-1}).
\end{array}$$
Thus, introducing the probability measure $\nu_q^N=\Phi_{q }\left(\eta_{q-1}^N\right)$ and the bounded function $\varphi=Q_{q+1,p}(f)$, our objective is to show that 
$$\sqrt{N}\ \nu_q^N\left(\left[Q_{q+1,\eta^N_q}-Q_{q+1}\right](\varphi)\right)\xrightarrow[N\to\infty]{\P}0.$$
Before going further, let us recall that if $G=\un_{S(\cdot)\geq L}$ is a potential function, $K$ a transition kernel and $M$ its truncated version defined by
$$
M(x,dy)=K(x,dy)~G(y)~+~K(1-G)(x)~\delta_x(dy),
$$ 
then for any finite measure $\mu$ and any bounded and measurable function $\varphi$, we have the following general formula
\begin{align}
\mu (G M(\varphi))=&\iint \mu(dy)G(y)K(y,x)G(x)\varphi(x)dx\nonumber\\
&+\iint \mu(dx)G(x)K(x,y)(1-G(y))\varphi(x)dy.\nonumber\\
=&\ \mu(G\times K[G\varphi])+\mu( K[1-G]\times (G\varphi)).
\end{align}
Thus, we get 
\begin{align}
&\nu_q^N\left(\left[Q_{q+1,\eta^N_q}-Q_{q+1}\right](\varphi)\right)\nonumber\\
&=\nu_q^N\left(G_{\eta_q^N}K_{q+1}[G_{\eta_q^N}\varphi]\right)- \nu_q^N\left(G_{q}K_{q+1}[G_{q}\varphi]\right)\nonumber\\
&\ \ + \nu_q^N\left(K_{q+1}[1-G_{\eta_q^N}]~(G_{\eta_q^N}\varphi)\right)-\nu_q^N\left(K_{q+1}[1-G_{q}](G_{q}\varphi)\right).\label{oamalkx}
\end{align} 
We may simplify a bit the latter by noticing that
$$\varphi=Q_{q+1,p}(f)=G_{q+1}\times \widetilde{Q}_{q+1,p}(f)=G_{\e_{q+1}}\times \widetilde{Q}_{q+1,p}(f).$$
Indeed, we know from Theorem \ref{kjsc} that 
$$L_{\eta_{q}^N}\xrightarrow[N\to\infty]{a.s.}L_{\eta_{q}}=L_q<L_{q+1}=L_{\e_{q+1}}.$$
Therefore, almost surely for $N>N_0$, we have $G_{\eta_q^N}\varphi=G_{q}\varphi=\varphi$, and (\ref{oamalkx}) reduces to
\begin{align}\label{blablabla}
\nu_q^N([Q_{q+1,\eta^N_q}-Q_{q+1}](\varphi))&=\nu_q^N((G_{\eta_q^N}-G_{q})K_{q+1}[\varphi])-\nu_q^N(K_{q+1}[G_{\eta_q^N}-G_{q}]\varphi).
\end{align}
In the remainder of the proof, we will only treat the more difficult case where the kernels $K_p$ are obtained by the Metropolis-Hastings procedure (see Section \ref{mh}) and we will suppose that Assumption [${\cal H}^a$] is satisfied. According to equation (\ref{kpaa}), we have
$$K_p[\varphi](x)=\ck_p^a[\varphi](x)+\ccr_p^a(x)\times\varphi(x).$$
In this expression, recall that $\ck_p^a$ has density
$$\ck_p^a(x,x')=a_p(x,x')\ck_p(x,x').$$
All the upcoming arguments remain valid in the easier case where $K_p$ itself has a density since it suffices to take $a=1$, so that $\ccr_p^a=0$ and $K_p=\ck_p$.\medskip

In the Metropolis-Hastings situation, combining (\ref{blablabla}) and (\ref{kpaa}), we are led to
\begin{align}
\nu_q^N([Q_{q+1,\eta^N_q}-Q_{q+1}](\varphi))&=\nu_q^N((G_{\eta_q^N}-G_{q})\ck^a_{q+1}[\varphi])-\nu_q^N(\ck^a_{q+1}[G_{\eta_q^N}-G_{q}]\varphi)\nonumber\\
&=A_q^N-B_q^N.\label{oamalkx-2}
\end{align}

Thanks to the coarea formula, $B_q^N$ rewrites
\begin{align}
B_q^N&=\ \int\nu_q^N(dx')\varphi(x')\int_{L_{q}}^{L_{q}^N}\left(\int_{S(x)=\ell}\ck_{q+1}^a(x',x)\frac{\dbar x}{|DS(x)|}\right)d\ell\nonumber\\
&=\ \nu_q^N\left(\varphi\int_{L_{q}}^{L_{q}^N}H_q^{1,a}(\cdot,\ell)d\ell\right).\nonumber
\end{align}
Next, since 
$$\nu_q^N=\Phi_{q }\left(\eta_{q-1}^N\right)=\frac{N}{\lceil N\alpha\rceil}\Pi_{q }(\eta_{q-1}^N),$$ 
we deduce from Assumption [${\cal H}^a$], Theorem \ref{kjsc}, Proposition \ref{asc} and Lemma \ref{undeplus} that 
\begin{align}
B_q^N&=\ (L_{q}^N-L_{q})\times\eta_q\left(\varphi H_q^{1,a}(\cdot,L_q)\right)+o_p(L_{q}^N-L_{q}),\nonumber\\
&=\ (L_{q}^N-L_{q})\iint_{S(x)=L_{q}}\e_q(dx')\varphi(x')\ck_{q+1}^a(x',x)\frac{\dbar x}{|DS(x)|}+o_p(L_{q}^N-L_{q}).\nonumber
\end{align}
Concerning $A_q^N$, coming back to (\ref{oamalkx-2}) and decomposing $\nu_q^N$ in absolutely continuous and discrete parts, we may write
\begin{align}
A_q^N&=\nu_q^{N,(0)}((G_{\eta_q^N}-G_{q})\ck^a_{q+1}[\varphi])+\nu_q^{N,(1)}((G_{\eta_q^N}-G_{q})\ck^a_{q+1}[\varphi])\nonumber\\
&=A_q^{N,(0)}+A_q^{N,(1)},\nonumber
\end{align}
where
\begin{equation}\label{ac}
\nu_q^{N,(0)}(dx)=\frac{1}{\lceil N\alpha\rceil}~\sum_{i=1}^{\lceil N\alpha\rceil}\ck^a_q(\widetilde{X}^i_{q-1},dx)G_{\eta_{q-1}^N}(x),
\end{equation}
and
\begin{equation}\label{d}
\nu_q^{N,(1)}(dx)=\frac{1}{\lceil N\alpha\rceil}~\sum_{i=1}^{\lceil N\alpha\rceil}(\ck^a_q[1-G_{\eta_{q-1}^N}](\widetilde{X}^i_{q-1})+\ccr_q^a(\widetilde{X}^i_{q-1}))\delta_{\widetilde{X}^i_{q-1}}(dx).
\end{equation}
As previously, since almost surely for $N>N_0$,
$$G_{\eta_{q-1}^N}(x)(G_{\eta_q^N}(x)-G_{q}(x))=G_{\eta_q^N}(x)-G_{q}(x),$$
we get
\begin{align}
A_q^{N,(0)}&=\int\frac{1}{\lceil N\alpha\rceil}~\sum_{i=1}^{\lceil N\alpha\rceil}\ck^a_q(\widetilde{X}^i_{q-1},dx)(\ck^a_{q+1}[\varphi]G_{\eta_{q-1}^N}(G_{\eta_q^N}-G_{q}))(x)\nonumber\\
&=\int\frac{1}{\lceil N\alpha\rceil}~\sum_{i=1}^{\lceil N\alpha\rceil}\ck^a_q(\widetilde{X}^i_{q-1},dx)(\ck^a_{q+1}[\varphi](G_{\eta_q^N}-G_{q}))(x)\nonumber\\
&=\ \tilde{\eta}_{q-1}^N\left(\int_{L_{q}}^{L_{q}^N}H_q^{\ck^a_{q+1}[\varphi],a}(\cdot,\ell)d\ell\right),\label{oamalkx-3bis}
\end{align}
the last equation consisting in the application of the coarea formula. Then, Assumption [${\cal H}$], Corollary \ref{tilde} and Lemma \ref{undeplus} yield
\begin{align}
A_q^{N,(0)}=&\ (L_{q}^N-L_{q})\iint_{S(x)=L_{q}}\e_{q}(dx')\ck^a_q(x',x)\ck^a_{q+1}[\varphi](x)\frac{\dbar x}{|DS(x)|}\nonumber\\
&+o_p(L_{q}^N-L_{q}).\nonumber
\end{align}
Using equation (\ref{kpa}), it is clear that for any pair $(x,x')$,
\begin{align}
\eta_{q}(dx')~\ck^a_q(x',x)\un_{S(x)\geq L_{q-1}}=&\alpha^{-q}~\eta(x')\un_{S(x')\geq L_{q-1}}dx'\ck^a_q(x',x)\un_{S(x)\geq L_{q-1}}\nonumber\\
=&\alpha^{-q}~\eta(x)\un_{S(x)\geq L_{q-1}}\ck^a_q(x,dx')\un_{S(x')\geq L_{q-1}}\nonumber.
\end{align}
Accordingly, denoting $w_{q-1}=\ck^a_q[1-G_{q-1}]$, this leads to
\begin{align}
A_q^{N,(0)}=&\ (L_{q}^N-L_{q})\int_{S(x)=L_{q}}(1-\ccr_{q}^a(x)-w_{q-1}(x))\alpha^{-q}\eta(x)\ck^a_{q+1}[\varphi](x)\frac{\dbar x}{|DS(x)|}\nonumber\\
&+o_p(L_{q}^N-L_{q}).\label{a}
\end{align}
By applying again (\ref{kpa}), and taking into account that $\varphi(x')\un_{S(x')\geq L_{q-1}}=\varphi(x')$, we have
$$\alpha^{-q}\ \eta(x)\ck^a_{q+1}[\varphi](x)=\int\e_{q}(dx')\ck^a_{q+1}(x',x)\varphi(x'),$$
and finally
\begin{align}
&A_q^{N,(0)}\nonumber\\
&=\ (L_{q}^N-L_{q})\iint_{S(x)=L_{q}}\e_{q}(dx')\varphi(x')\ck^a_{q+1}(x',x)(1-\ccr_{q}^a(x)-w_{q-1}(x))\frac{\dbar x}{|DS(x)|}\nonumber\\
&\ \ +o_p(L_{q}^N-L_{q}).\nonumber
\end{align}
Next, we come back to $A_q^{N,(1)}$, defined as
\begin{align}
A_q^{N,(1)}&=\frac{1}{\lceil N\alpha\rceil}\sum_{i=1}^{\lceil N\alpha\rceil}\left((\ck^a_q[1-G_{\eta_{q-1}^N}]+\ccr_q^a)(\ck^a_{q+1}[\varphi](G_{\eta_q^N}-G_{q}))\right)(\widetilde{X}^i_{q-1})\nonumber\\
&=\tilde\eta_{q-1}^N\left((\ck^a_q[1-G_{\eta_{q-1}^N}]+\ccr_q^a)\ck^a_{q+1}[\varphi](G_{\eta_q^N}-G_{q})\right).\label{oamalkxx}
\end{align}
Then, if we denote 
$$w_{q-1}^N(x)=\ck^a_q[1-G_{\eta_{q-1}^N}](x)=1-\ccr_q^a(x)-\ck^a_q [G_{\eta_{q-1}^N}](x),$$
we have
\begin{equation}\nonumber
A_q^{N,(1)}=\tilde\eta_{q-1}^N\left((w_{q-1}^N+\ccr_q^a)\ck^a_{q+1}[\varphi](G_{\eta_q^N}-G_{q})\right).
\end{equation}
At this step, it is quite natural to consider the deterministic functions $w_{q-1}^{\delta^-}\leq w_{q-1}^{\delta^+}$ defined by
$$w_{q-1}^{\delta^\pm}(x)=\ck^a_q[1-G_{L_{q-1}\pm\delta}](x)=\ck^a_q[1-G_{L_{q-1}\pm\delta}](x).$$
Accordingly, let us also introduce the random variable
\begin{equation}\label{A}
\hat A_q^{N,(1)}=\tilde\eta_{q-1}^N\left((w_{q-1}+\ccr_q^a)\ck^a_{q+1}[\varphi](G_{\eta_q^N}-G_{q})\right).
\end{equation} 
In what follows, we assume that $f$ is non-negative, otherwise we decompose $f=f^+-f^-$ and the same reasoning applies to both parts. If $f\geq 0$, then the same is true for $\varphi=Q_{q+1,p}(f)$ and we have $0\leq \ck^a_{q+1}[\varphi]\leq 1$. Besides, we remark that the sign of $w_{q-1}^N(x)-w_{q-1}(x)$ is independent of $x$, which is also true for $G_{\eta_q^N}(x)-G_{q}(x)$. As a consequence, since $L_{q-1}^N$ tends almost surely to $L_{q-1}$, we have that, almost surely for $N>N_0$,
$$\left|A_q^{N,(1)}-\hat A_q^{N,(1)}\right|\leq\left|\Delta_{q-1}^N\right|,$$
where  
\begin{equation}\nonumber
\Delta_{q-1}^N=\tilde\eta_{q-1}^N\left((w_{q-1}^{\delta^+}-w_{q-1}^{\delta^-})(G_{\eta_q^N}-G_{q})\right).
\end{equation}
We will first focus our attention on $\Delta_{q-1}^N$ and then exhibit the limit of $\hat A_q^{N,(1)}$. Concerning $\Delta_{q-1}^N$, we may reformulate it as
$$\Delta_{q-1}^N=\frac{N}{\lceil N\alpha\rceil}~\e_{q-1}^N\left((w_{q-1}^{\delta^+}-w_{q-1}^{\delta^-})(G_{\eta_q^N}-G_{q})\right),$$
and Corollary \ref{equivalence} implies that
\begin{equation}\label{abc}
\Delta_{q-1}^N=\frac{1}{\alpha}~\nu_{q-1}^N\left((w_{q-1}^{\delta^+}-w_{q-1}^{\delta^-})(G_{\eta_q^N}-G_{q})\right)+o_p(1/\sqrt{N}).
\end{equation}
As before, given ${\cal F}_{q-2}^N$, we split 
$$\nu_{q-1}^N=\Phi_{q-1}(\eta_{q-2}^N)=\nu_{q-1}^{N,(0)}+\nu_{q-1}^{N,(1)}$$ 
in absolutely continuous and discrete parts, see equations (\ref{ac}) and (\ref{d}) with $(q-1)$ instead of $q$, $\ck^a_{q-1}$ instead of $\ck^a_{q}$ and $\ccr^a_{q-1}$ instead of $\ccr^a_{q}$, leading to 
\begin{equation}\nonumber
\Delta_{q-1}^{N}=\frac{1}{\alpha}~\left(\Delta_{q-1}^{N,(0)}+\Delta_{q-1}^{N,(1)}\right)+o_p(1/\sqrt{N}),
\end{equation}
where
$$\Delta_{q-1}^{N,(0)}=\int\frac{1}{\lceil N\alpha\rceil}\sum_{i=1}^{\lceil N\alpha\rceil}\ck^a_{q-1}(\widetilde{X}^i_{q-2},dx)((w_{q-1}^{\delta^+}-w_{q-1}^{\delta^-})(G_{\eta_q^N}-G_{q}))(x),$$
and
$$\Delta_{q-1}^{N,(1)}=\frac{1}{\lceil N\alpha\rceil}\sum_{i=1}^{\lceil N\alpha\rceil}((w_{q-2}^N+\ccr^a_{q-1})(w_{q-1}^{\delta^+}-w_{q-1}^{\delta^-})(G_{\eta_q^N}-G_{q}))(\widetilde{X}^i_{q-2}).$$
Clearly, $\Delta_{q-1}^{N,(0)}$ shares some resemblance with $A_q^{N,(0)}$ as given in (\ref{oamalkx-3bis}). Therefore, mutatis mutandis, we get an equivalent expression as (\ref{a}), namely
\begin{align}
&\Delta_{q-1}^{N,(0)}\nonumber\\
&=(L_{q}^N-L_{q})\int_{S(x)=L_{q}}\e_{q-1}(x)((w_{q-1}^{\delta^+}-w_{q-1}^{\delta^-})(1-\ccr^a_{q-1}-w_{q-2}))(x)\ \frac{\dbar x}{|DS(x)|}\nonumber\\
&\ \ +o_p(L_{q}^N-L_{q}).\nonumber
\end{align}
Since $0\leq 1-\ccr^a_{q-1}-w_{q-2}=\ck^a_{q-1}[G_{q-2}]\leq 1$, we deduce in particular that
\begin{align}
\left|\Delta_{q-1}^{N,(0)}\right|\leq&\left|L_{q}^N-L_{q}\right|\int_{S(x)=L_{q}}\e_{q-1}(x)(w_{q-1}^{\delta^+}-w_{q-1}^{\delta^-})(x)\ \frac{\dbar x}{|DS(x)|}\nonumber\\
&+o_p(L_{q}^N-L_{q}).\label{delta}
\end{align}
Regarding $\Delta_{q-1}^{N,(1)}$, since $0\leq w_{q-2}^N\leq 1$, we get $\left|\Delta_{q-1}^{N,(1)}\right|\leq\left|\Delta_{q-2}^{N}\right|$, with
$$\Delta_{q-2}^{N}=\frac{1}{\lceil N\alpha\rceil}\sum_{i=1}^{\lceil N\alpha\rceil}((w_{q-1}^{\delta^+}-w_{q-1}^{\delta^-})(G_{\eta_q^N}-G_{q}))(\widetilde{X}^i_{q-2}).$$
Putting all pieces together yields
$$\left|A_q^{N,(1)}-\hat A_q^{N,(1)}\right|\leq\frac{1}{\alpha}\left(\left|\Delta_{q-1}^{N,(0)}\right|+\left|\Delta_{q-2}^{N}\right|\right)+o_p(1/\sqrt{N}),$$
and finally
$$\left|A_q^{N,(1)}-\hat A_q^{N,(1)}\right|\leq\alpha^{-1}\left|\Delta_{q-1}^{N,(0)}\right|+\dots+\alpha^{1-q}\left|\Delta_{1}^{N,(0)}\right|+\alpha^{1-q}\left|\Delta_{0}^{N}\right|+o_p(1/\sqrt{N}).$$
By (\ref{delta}), for every $k\in\{1,\dots,q-1\}$, we have the upper-bound
\begin{align}
\left|\Delta_{k}^{N,(0)}\right|\leq&\left|L_{q}^N-L_{q}\right|\int_{S(x)=L_{q}}\e_{k}(x)(w_{q-1}^{\delta^+}-w_{q-1}^{\delta^-})(x)\ \frac{\dbar x}{|DS(x)|}\nonumber\\
&+o_p(L_{q}^N-L_{q}),\nonumber
\end{align}
and, by (\ref{abc}), we have
$$\Delta_{0}^{N}=\frac{1}{\alpha}~\nu_{0}^N\left((w_{q-1}^{\delta^+}-w_{q-1}^{\delta^-})(G_{\eta_q^N}-G_{q})\right)+o_p(1/\sqrt{N}).$$
Since $\nu_{0}^N=\eta_0=\eta$, the coarea formula yields
\begin{align}
\left|\Delta_{0}^{N}\right|\leq&\frac{1}{\alpha}\left|L_{q}^N-L_{q}\right|\int_{S(x)=L_{q}}\e_0(x)(w_{q-1}^{\delta^+}-w_{q-1}^{\delta^-})(x)\ \frac{\dbar x}{|DS(x)|}\nonumber\\
&+o_p(L_{q}^N-L_{q})+o_p(1/\sqrt{N}).\nonumber
\end{align}
Lebesgue's dominated convergence theorem ensures that 
$$\int_{S(x)=L_{q}}\e_{k}(x)(w_{q-1}^{\delta^+}-w_{q-1}^{\delta^-})(x)\ \frac{\dbar x}{|DS(x)|}\xrightarrow[\delta\to 0]{}0,$$
and Lemma \ref{L} says that $L_{q}^N-L_{q}={\cal O}_p(1/\sqrt{N})$, so we conclude that 
$$A_q^{N,(1)}-\hat A_q^{N,(1)}=o_p(1/\sqrt{N}).$$
Now we turn to the estimation of $\hat A_q^{N,(1)}$ as defined in (\ref{A}). The analysis is roughly the same as for $\Delta_{q-1}^N$ except that we have to be a bit more precise since this time we want an estimate and not an upper-bound. However, we can reformulate it as
$$\hat A_q^{N,(1)}=\frac{N}{\lceil N\alpha\rceil}~\e_{q-1}^N\left((\ccr^a_q+w_{q-1})\ck^a_{q+1}[\varphi](G_{\eta_q^N}-G_{q})\right),$$
and Corollary \ref{equivalence} implies that
$$\hat A_q^{N,(1)}=\frac{1}{\alpha}~\nu_{q-1}^N\left((\ccr^a_q+w_{q-1})\ck^a_{q+1}[\varphi](G_{\eta_q^N}-G_{q})\right)+o_p(1/\sqrt{N}).$$
Again, given ${\cal F}_{q-2}^N$, we split $\nu_{q-1}^N=\nu_{q-1}^{N,(0)}+\nu_{q-1}^{N,(1)}$ into its absolutely continuous and discrete parts to get
$$\hat A_q^{N,(1)}=\frac{1}{\alpha}\left(A_{q-1}^{N,(0)}+A_{q-1}^{N,(1)}\right)+o_p(1/\sqrt{N}),$$
where, as in (\ref{oamalkx-3bis}) and (\ref{oamalkxx}),
$$A_{q-1}^{N,(0)}=\int\frac{1}{\lceil N\alpha\rceil}~\sum_{i=1}^{\lceil N\alpha\rceil}\ck^a_{q-1}(\widetilde{X}^i_{q-2},dx)((\ccr^a_q+w_{q-1})\ck^a_{q+1}[\varphi](G_{\eta_q^N}-G_{q}))(x),$$
and
$$A_{q-1}^{N,(1)}=\frac{1}{\lceil N\alpha\rceil}\sum_{i=1}^{\lceil N\alpha\rceil}(\ck^a_{q-1}(1-G_{\eta_{q-2}^N})+\ccr^a_{q-1})((\ccr^a_q+w_{q-1})\ck^a_{q+1}[\varphi](G_{\eta_q^N}-G_{q}))(\widetilde{X}^i_{q-2})$$
By the same arguments as above, under Assumption [${\cal H}^a$], it is readily seen that
\begin{align}
&A_{q-1}^{N,(0)}=\ (L_{q}^N-L_{q})\nonumber\\
&\iint_{S(x)=L_{q}}\e_{q-1}(dx')\varphi(x')\ck^a_{q+1}(x',x)((\ccr^a_q+w_{q-1})(1-\ccr^a_{q-1}-w_{q-2}))(x)\frac{\dbar x}{|DS(x)|}\nonumber\\
&\ +o_p(L_{q}^N-L_{q}).\nonumber
\end{align}
Moreover, by the same machinery as for the majorization of $\Delta_{q-1}^N$, we get
$$A_{q-1}^{N,(1)}-\hat A_{q-1}^{N,(1)}=o_p(1/\sqrt{N}).$$
Consequently, we have 
$$A_q^N=A_q^{N,(0)}+\frac{1}{\alpha}A_{q-1}^{N,(0)}+\frac{1}{\alpha}\hat A_{q-1}^{N,(1)}+o_p(1/\sqrt{N}).$$
At this point, it remains to notice that
$$\eta_q(dx')\varphi(x')=\frac{1}{\alpha}\eta_{q-1}(dx')\varphi(x'),$$
which implies that
\begin{align}
&A_{q}^{N,(0)}+\frac{1}{\alpha}A_{q-1}^{N,(0)}=(L_{q}^N-L_{q})\nonumber\\
&\iint_{S(x)=L_{q}}\e_{q}(dx')\ck^a_{q+1}(x',x)((1-\ccr^a_q-w_{q-1})(\ccr^a_{q-1}+w_{q-2}))(x)\frac{\dbar x}{|DS(x)|}\nonumber\\
&\ \ +o_p(L_{q}^N-L_{q}),\nonumber
\end{align} 
and a straightforward recursion gives
\begin{align}
&A_{q}^{N}=\ (L_{q}^N-L_{q})\nonumber\\
&\iint_{S(x)=L_{q}}\e_{q}(dx')\varphi(x')\ck^a_{q+1}(x',x)((1-\ccr^a_q-w_{q-1})\dots(\ccr^a_1+ w_0))(x)\frac{\dbar x}{|DS(x)|}\nonumber\\
&+\alpha^{1-q}\hat A_{1}^{N,(1)}+o_p(L_{q}^N-L_{q})+o_p(1/\sqrt{N}),\nonumber
\end{align}
where
$$\hat A_1^{N,(1)}=\frac{1}{\alpha^q}~\nu_{0}^N\left((\ccr^a_q+w_{q-1})\dots (\ccr^a_1+w_0)\ck_{q+1}^a[\varphi](G_{\eta_q^N}-G_{q})\right)+o_p(1/\sqrt{N}).$$
Since $\nu_{0}^N=\eta$, we finally get
\begin{align}
&\hat A_1^{N,(1)}=\ (L_{q}^N-L_{q})\nonumber\\
&\iint_{S(x)=L_{q}}\e_{q}(dx')\varphi(x')\ck^a_{q+1}(x',x)((\ccr^a_q+w_{q-1})\dots (\ccr^a_1+w_0))(x)\frac{\dbar x}{|DS(x)|}\nonumber\\
&+o_p(L_{q}^N-L_{q})+o_p(1/\sqrt{N}),\nonumber
\end{align}
so that, coming back to (\ref{oamalkx-2}) and thanks to Proposition \ref{L}, we have eventually shown that
$$\nu_q^N([Q_{q+1,\eta^N_q}-Q_{q+1}](\varphi))=o_p(L_{q}^N-L_{q})+o_p(1/\sqrt{N})=o_p(1/\sqrt{N}).$$
This terminates the proof of Proposition \ref{VieillePuteAlbanaise}.
\hfill  $\blacksquare$
\medskip

The following lemma is a key tool to prove Proposition \ref{L} and its Corollary \ref{equivalence}, which were useful in the previous proof.

\begin{lem}\label{techlemma}
For any $C>0$, for any integer $0\leq q<n$ and for any $L\in\{L_{q},\dots,L_{n-1}\}$, consider the class of sets 
$$\A_{N,C}=\left\{S^{-1}\left(\left[L-\frac{c_1}{\sqrt{N}},L+\frac{c_2}{\sqrt{N}}\right]\right), 0<c_1<C, 0<c_2<C\right\}.$$
Then, for any $\phi\in\Ba(\R^d)$, we have that 
$$\sup_{A\in\A_{N,C}}\sqrt{N}\left|\nu_q^N(\phi\1_A)-\eta_q^N(\phi\1_A)\right|\xrightarrow[N\to\infty]{\P} 0.$$ 
\end{lem}

\paragraph{Proof}
Here again, the proof is made given ${\cal F}_{q-1}^N$. Let $A_{N,C}$ denote the largest set in $\A_{N,C}$, i.e. 
$$A_{N,C}=S^{-1}\left(\left[L-\frac{C}{\sqrt{N}},L+\frac{C}{\sqrt{N}}\right]\right).$$ 
Let us write some preliminary algebra. In the following, $k_N$ stands for the number of sample points belonging to $A_{N,C}$, meaning that
$$k_N=N\times\eta_q^N(A_{N,C})=\sum_{i=1}^N \1_{A_{N,C}}(X_q^i).$$
We start from the decomposition
\begin{eqnarray}
\lefteqn{\nonumber
\sup_{A\in\A_{N,C}}\sqrt{N}\left|\nu_q^N(\phi\1_A)-\eta_q^N(\phi\1_A)\right|}\\
&\leq&\sqrt{N}\ \nu_q^N(A_{N,C})\sup_{A\in\A_{N,C}} \left|\frac{\nu_q^N(\phi\1_A)}{\nu_q^N(A_{N,C})}-\frac{1}{k_N}\sum_{i=1}^N \1_A(X^i_q)\phi(X^i_q)\right| \label{lemsup1}\\
&&+\sup_{A\in\A_{N,C}} \left|\frac{1}{k_N}\sum_{i=1}^N \1_A(X^i_q)\phi(X^i_q)\right|\times\left|\sqrt{N}\ \nu_q^N(A_{N,C})-\frac{k_N}{\sqrt{N}}\right|.\label{lemsup2}
\end{eqnarray}

\medskip

Consider first expression~(\ref{lemsup1}). We study the class $\A_{N,C}$ from the viewpoint of Vapnik-Chervonenkis theory (see for example Chapters 12 and 13 in \cite{DGL}). We denote by $s(\A_{N,C},N)$ the shattering coefficient of $\A_{N,C}$. Very elementary reasoning gives that $s(\A_{N,C},N)\leq N^2$.\medskip

As $\phi$ is bounded, for any $\varepsilon>0$ we can find a simple function $\phi^\varepsilon=\sum_{j=1}^{n_\varepsilon} b_j \1_{B_j}$ such that $\|\phi-\phi^\varepsilon\|<\varepsilon$. Let us denote by $\B_\varepsilon$ the finite collection of Borelian sets in the expression of $\phi^\varepsilon$. If we consider now 
$$\A_{N,C}^\varepsilon=\left\{A=A_1\cap A_2, A_1\in\A_{N,C}, A_2\in \B_\varepsilon\right\},$$
then it is clear that its shatter coefficient verifies $s(\A_{N,C}^\varepsilon,N)\leq 2^{n_\varepsilon} N^2$.\medskip 

Now, in (\ref{lemsup1}), we show that the supremum factor goes to $0$ in probability. We first have
\begin{align}
&\sup_{A\in\A_{N,C}} \left|\frac{\nu_q^N(\phi\1_A)}{\nu_q^N(A_{N,C})}-\frac{1}{k_N}\sum_{i=1}^N \1_A(X^i_q)\phi(X^i_q)\right|\nonumber\\
&\ \ \leq\sup_{A\in\A_{N,C}}\left|\frac{\nu_q^N((\phi-\phi^\varepsilon)\1_A)}{\nu_q^N(A_{N,C})}\right|+\sup_{A\in\A_{N,C}} 
\left|\frac{\nu_q^N(\phi^\varepsilon\1_A)}{\nu_q^N(A_{N,C})}-\frac{1}{k_N}\sum_{i=1}^N \1_A(X^i_q)\phi^\varepsilon(X^i_q)\right|\nonumber\\
&\ \ \ \ +\sup_{A\in\A_{N,C}} 
\left|\frac{1}{k_N}\sum_{i=1}^N \1_A(X^i_q)(\phi^\varepsilon-\phi)(X^i_q)\right|,\nonumber
\end{align}
hence
\begin{align}
&\sup_{A\in\A_{N,C}} \left|\frac{\nu_q^N(\phi\1_A)}{\nu_q^N(A_{N,C})}-\frac{1}{k_N}\sum_{i=1}^N \1_A(X^i_q)\phi(X^i_q)\right|\nonumber\\
&\leq 2 \|\phi-\phi^\varepsilon\| + \sup_{A\in\A_{N,C}}  \left|\sum_{j=1}^{n_\varepsilon} b_j\left(\frac{\nu_q^N(\1_{A\cap B_j})}{\nu_q^N(A_{N,C})} - \frac{1}{k_N}\sum_{i=1}^N \1_{A\cap B_j}(X^i_q)\right)\right|\nonumber\\
&\leq 2\varepsilon + \left(\sum_{j=1}^{n_\varepsilon} |b_j|\right) \times \sup_{A\in\A_{N,C}^{\varepsilon}} \left|\frac{\nu_q^N(\1_A)}{\nu_q^N(A_{N,C})}-\frac{1}{k_N}\sum_{i=1}^N \1_A(X^i_q)\right|\nonumber\\
&\leq 2\varepsilon + \left(\sum_{j=1}^{n_\varepsilon} |b_j|\right) \times \varepsilon^\prime\nonumber\\
&\leq 3\varepsilon,\nonumber
\end{align}
for $\varepsilon^\prime$ chosen small enough, with probability at least 
$$
1-8s(\A_{N,C}^\varepsilon,N) e^{-N\varepsilon^{\prime 2}/32}\geq 1-2^{n_\varepsilon+3} N^2 e^{-N\varepsilon^{\prime 2 }/32},
$$
which can be made arbitrarily close to $1$  for $N$ large enough. We notice that here we have used Theorem 12.5 in \cite{DGL}, and the fact that, given ${\cal F}_{q-1}^N$, the $X_q^i$'s are i.i.d. with distribution $\nu_q^N$, and thus the $k_N$  ones in $A_{N,C}$ are i.i.d. with distribution $\nu_q^N.\1_{A_{N,C}}/\nu_q^N(A_{N,C})$. 
\medskip

Now, to complete the proof of the lemma, it suffices to show that the pre-factor $\sqrt{N}\nu_q^N(A_{N,C})$ in (\ref{lemsup1}) can be bounded with arbitrarily large probability. In this aim, we proceed by induction on $q$. Consider first $q=0$. In that case $\nu_q^N=\eta$, and it is clear using the coarea formula and the law of large numbers that 
$$\nu_q^N(A_{N,C})={\cal O}_p(1/\sqrt{N}).$$
For the general case $q>0$, we have the decomposition $\nu_q^N=\nu_q^{N,(0)} + \nu_q^{N,(1)}$ where the first term is absolutely continuous with respect to Lebesgue's measure, and the second term is a discrete one. A quick inspection reveals that
\begin{equation}\label{nudec}
\nu_q^{N,(0)}\leq \frac{1}{\alpha}\eta_{q-1}^N \ck^a_q\hspace{1cm}\mbox{and}\hspace{1cm}\nu_q^{N,(1)}\leq \frac{1}{\alpha} \eta_{q-1}^N.
\end{equation}
When applied to $A_{N,C}$ both are ${\cal O}_p(1/\sqrt{N})$. For the first one we simply apply the coarea formula and the law of large numbers. For the second one, we notice that $N \eta_{q-1}^N(A_{N,C})$ is a Binomial r.v. with parameters $N$ and $\nu_{q-1}^N(A_{N,C})$. The mean $\nu_{q-1}^N(A_{N,C})$ is ${\cal O}_p(1/\sqrt{N})$ by the induction assumption. For the distance to the mean we use Hoeffding's inequality
$$
\P\left(\left.\left| \nu_{q-1}^N(A_{N,C}) - \eta_{q-1}^N(A_{N,C}) \right| \geq \frac{A}{\sqrt{N}}\right|{\cal F}_{q-2}^N \right)\leq 2e^{-2A^2},
$$
which can be made arbitrarily small by choosing $A$ large enough. This shows that $| \nu_{q-1}^N(A_{N,C}) - \eta_{q-1}^N(A_{N,C}) |$ is also ${\cal O}_p(1/\sqrt{N})$.\medskip

Consider now expression~(\ref{lemsup2}). It is clear that the supremum is less than $\|\phi\|$. For the factor $|\sqrt{N}\nu_q^N(A_{N,C})-\frac{k_N}{\sqrt{N}}|$, let us denote $I_q^N=\nu_q^N(A_{N,C})$. From usual considerations on the $X^i_q$'s, we see that ${k_N}$ is Binomial ${\cal B}(N,I_q^N)$ distributed, thus  we have
$$\E\left[\left.\frac{k_N}{N}\right|{\cal F}_{q-1}^N\right]=I_q^N\hspace{1cm}\mbox{and}\hspace{1cm}\V\left(\left.\frac{k_N}{N}\right|{\cal F}_{q-1}^N\right)=\frac{I_q^N(1-I_q^N)}{N}.$$ 
By Chebyshev's inequality we deduce that, for any $\varepsilon>0$, 
$$\P\left(\left.\sqrt{N}\left|\frac{k_N}{N}-I_q^N\right|>\varepsilon\right|{\cal F}_{q-1}^N\right)\leq \frac{I_q^N(1-I_q^N)}{\varepsilon^2}\xrightarrow[N\to\infty]{\P}0,$$
since, as justified above, $I_q^N=\nu_q^N(A_{N,C})={\cal O}_P(1/\sqrt{N})$.
\hfill  $\blacksquare$

\begin{pro}\label{L}
For all $q\in\{0,\dots,n-1\}$,
$$L_q^N-L_{q} = {\cal O}_p(1/\sqrt{N}).$$
\end{pro}

\paragraph{Proof}
The proof is done by induction on $q$. We will actually make the induction on the following double property: for all $\delta>0$, for all measurable function $\phi$ such that $0\leq\phi\leq 1$ and with support above $L_q$ (i.e. $\phi=G_{q} \phi$), there exist $C>0$ and $N_0$ such that for all $N>N_0$, with probability at least $(1-\delta)$, we have
$$\left|L_q^N-L_{q}\right|\leq \frac{C}{\sqrt{N}}\hspace{1cm}\mbox{and}\hspace{1cm}\left|\left(\eta_{q}-\nu_{q}^N\right)(\phi)\right|\leq \frac{C}{\sqrt{N}}.$$
First note that for $q=0$, since $\nu_{0}^N=\eta_{0}$, the second assertion is trivial, and the first one is obtained by very standard properties of empirical quantiles (e.g. CLT) when the i.i.d. sample is drawn from a distribution with a strictly positive density at point $L_0$. 
\medskip
 
Now, assume the property is true up to step $(q-1)$. Then, by we have
\begin{align}
\alpha\left(\nu_q^N-\eta_q\right)(\phi)
=&\ \eta_{q-1}^N\left(G_{L_{q-1}^N}M_{q,\eta_{q-1}^N}\phi\right) - \nu_{q-1}^N\left(G_{q-1}M_{q}\phi\right)\label{tl1}\\
&+  \left(\nu_{q-1}^N- \eta_{q-1}\right)(G_{q-1}M_{q}\phi)+o(1/\sqrt{N}).\label{tl2}
\end{align}
The second term (\ref{tl2}) is easy as $\|G_{q-1}M_{q}\phi\|\leq 1$ and, from the recurrence assumption, its absolute value is less than $C/\sqrt{N}$ with probability at least $(1-\delta)$.
\medskip

For the first term, namely (\ref{tl1}), let us write
\begin{eqnarray}
\lefteqn{\nonumber
\left|\eta_{q-1}^N\left(G_{L_{q-1}^N}M_{q,\eta_{q-1}^N}\phi\right) - \nu_{q-1}^N\left(G_{q-1}M_q\phi\right) \right|}\\
\label{tl3}
&\leq&
\left|\eta_{q-1}^N\left(\left(G_{L_{q-1}^N}M_{q,\eta_{q-1}^N}-G_{q-1}M_q\right)\phi\right)\right|\\
\label{tl4}
&&+
\left|\nu_{q-1}^N\left(G_{q-1}M_q\phi\right)-\eta_{q-1}^N \left(G_{q-1}M_q\phi\right)\right|.
\end{eqnarray}
Let us first consider (\ref{tl4}). Since $\eta_{q-1}^N$ is an empirical measure of an i.i.d. sample drawn with $\nu_{q-1}^N$, Chebyshev's inequality implies that, for all $t>0$,
$$\P\left(\left|\nu_{q-1}^N\left(G_{q-1}M_q\phi\right)-\eta_{q-1}^N \left(G_{q-1}M_q\phi\right)\right|\geq t\sigma_N\right)\leq \frac{1}{t^2},$$
with 
$$\sigma_N\leq\frac{1}{\sqrt{N}}\sqrt{\nu_{q-1}^N\left[(G_{q-1}M_q\phi)^2\right]}\leq\frac{1}{\sqrt{N}}.$$
Thus, if we take 
$$t=1/\sqrt{\delta}\hspace{1cm}\mbox{and}\hspace{1cm}C>\frac{1}{\sqrt{\delta}},$$
it turns out that, for $N$ large enough, we have with probability at least $(1-\delta)$, 
$$\left|\nu_{q-1}^N\left(G_{q-1}M_q\phi\right)-\eta_{q-1}^N\left(G_{q-1}M_q\phi\right)\right|\leq\frac{C}{\sqrt{N}}.$$
Now we decompose (\ref{tl3}) in a similar way as (\ref{oamalkx-2}) and taking into account that $G_q\phi=\phi$, which gives
\begin{align}
&\eta_{q-1}^N\left(\left(G_{L_{q-1}^N}M_{q,L_{q-1}^N}-G_{q-1}M_{q}\right)\phi\right)\nonumber\\
&=\eta_{q-1}^N \left(\left(G_{L_{q-1}^N}-G_{q-1}\right)\ck^a_q[\phi]\right) 
- \eta_{q-1}^N\left(\phi \ \ck^a_q[G_{L_{q-1}^N}-G_{q-1}]\right).\label{blabla}
\end{align}
With probability at least $(1-\delta)$, for $N$ large enough, we have for the second term, using the recurrence assumption and the coarea formula,
\begin{eqnarray*}
\lefteqn{
\left|\eta_{q-1}^N\left(\phi \ \ck^a_q[G_{L_{q-1}^N}-G_{q-1}]\right)\right|}\\
&\leq&
\left|\eta_{q-1}^N\left(\phi \ \ck^a_q[G_{L_{q-1}-\frac{C}{\sqrt{N}} }   -G_{L_{q-1}+\frac{C}{\sqrt{N}}}]\right)\right|\\
&\leq& \left|\eta_{q-1}^N\left(\phi \int_{\{S(y)=L_{q-1}\}} \ck_q^a(\cdot,y)\frac{\dbar y}{|DS(y)|}\right)\right|\times\frac{2C}{\sqrt{N}} + o_p(1/\sqrt{N}),
\end{eqnarray*}
with the main factor converging in probability to 
$$\eta_{q-1}\left(\phi \int_{\{S(y)=L_{q-1}\}} \ck_q^a(\cdot,y)\frac{\dbar y}{|DS(y)|}\right).$$
For the first term in (\ref{blabla}) we have, thanks to Lemma \ref{techlemma}, 
$$\eta_{q-1}^N \left(\left(G_{L_{q-1}^N}-G_{q-1}\right)\ck_q^a[\phi]\right) 
= \nu_{q-1}^N \left(\left(G_{L_{q-1}^N}-G_{q-1}\right)\ck_q^a[\phi]\right) + o_p(1/\sqrt{N}).$$
We then upper-bound $\nu_{q-1}^N$ like in (\ref{nudec}) in order to write
\begin{align}
&\left|\nu_{q-1}^N \left(\left(G_{L_{q-1}^N}-G_{q-1}\right) \ck^a_q[\phi]\right)\right|\nonumber\\
&\leq\frac{1}{\alpha}\left|\eta_{q-2}^N \left(\left(G_{L_{q-1}^N}-G_{q-1}\right) \ck^a_q[\phi]\right)\right|\nonumber \\
&\ \ +\frac{1}{\alpha}\left|\eta_{q-2}^N\left(\int \ck^a_{q-1}(\cdot,y)(G_{L_{q-1}^N}-G_{q-1})(y)\ck^a_q[\phi](y) dy\right)\right|.\nonumber
\end{align}
For the second term, we use the coarea formula and the recurrence assumption just as above, and for the first term, we replace $\eta_{q-2}^N$ with $\nu_{q-2}^N$ by virtue of Lemma \ref{techlemma}. We iterate the reasoning until we get terms with $\nu_{0}^N=\eta$, which can be dealt by applying the coarea formula again. 
\medskip

Now we consider the other part of the recurrence assumption. Let us define the function  $F_N(\ell)=1-\nu_q^N(G_\ell)$ and 
$$L_{\nu_q^N}=\inf \{t\ \mbox{such that}\ (1-F_N(t))\geq 1-\alpha\}.$$ 
Following the same arguments as in the proof of Theorem \ref{kjsc}, we can easily see that $L_{\nu_q^N}$ a.s. converges to $L_q$. We obviously have
\begin{equation}\label{truc}
\left|L_q^N-L_q\right|\leq \left|L_q^N - L_{\nu_q^N}\right| + \left|L_{\nu_q^N} - L_q\right|.
\end{equation}
We first deal with $|L_q^N - L_{\nu_q^N}|$. 
From the proof of Lemma \ref{lem.TCL.1} we see that 
$$F_N(L_{\nu_q^N})=\alpha+o_{L^2}(1/\sqrt{N}),$$ 
so that 
$$\V(F_N(L_q^N)-F_N(L_{\nu_q^N})|\ \F_{q-1}^N )\leq 2\ \E\left[\left.\left(  F_N(L_q^N) -\alpha\right)^2\ \right| \ \F_{q-1}^N \right] +o_p(1/N).$$
Moreover, from Lemma \ref{lem.TCL.1} we have
\begin{equation}\label{convergence}
N~\E\left[\left.\left(  F_N(L_q^N) -\alpha\right)^2\ \right| \ \F_{q-1}^N \right] \xrightarrow[N\to\infty]{\P} \alpha(1-\alpha),
\end{equation}
Hence, using Chebyshev's inequality we see that, given $\F_{q-1}^N$, the random variable $\sqrt{N}(F_N(L_q^N)-F_N(L_{\nu_q^N}))$ is bounded with arbitrarily large probability, and so it is unconditionally,  for in (\ref{convergence}) the limit is deterministic.\medskip

As mentioned before, the function $F_N$ is absolutely continuous except at a finite number of points, namely at most the $\lceil N\alpha \rceil$ largest $Y_{q-1}^i$'s. Denoting $f_N$ the density of the absolutely continuous part of $F_N$, and $J_i$'s the heights of the jumps, we may write
$$F_N(L_{\nu_q^N})-F_N(L_q^N)=\int_{L_q^N}^{L_{\nu_q^N}}f_N(\ell)d\ell+\sum_{i:Y_{q-1}^i\in[L_q^N,L_{\nu_q^N}]}J_i,$$
where $[L_q^N,L_{\nu_q^N}]$ stands for $[L_q^N,L_{\nu_q^N}]$ or $[L_{\nu_q^N},L_q^N]$. 
We want to show that, with large probability,
$$
|F_N(L_{\nu_q^N})-F_N(L_q^N)|\geq \left|\int_{L_q^N}^{L_{\nu_q^N}}f_N(\ell)d\ell\right|\geq C_q |L_q^N-L_{\nu_q^N}|,
$$
where $C_q>0$ is some deterministic constant.
We have
$$
\int_{L_q^N}^{L_{\nu_q^N}}f_N(\ell)d\ell = \tilde\eta_{q-1}^N\int_{L_q^N}^{L_{\nu_q^N}} H^{1,a}_q(.,\ell)\ d\ell,
$$
Therefore, using assumption $[{\cal H}^a]$, as for $N$ large both $L_q^N$ and $L_{\nu_q^N}$ are close to $L_q$,  we can write
$$
-\varepsilon h(.)\leq H^{1,a}_q(.,\ell) - H^{1,a}_q(.,L_q)\leq \varepsilon h(.),
$$
uniformly for $\ell$ between $L_q^N$ and $L_{\nu_q^N}$.
From all that we get
$$
|F_N(L_{\nu_q^N})-F_N(L_q^N)|\geq  |L_q^N-L_{\nu_q^N}| \times \tilde\eta_{q-1}^N(H^{1,a}_q(.,L_q) - 2\varepsilon h).
$$
By the law of large numbers, the last factor on the right can be made larger than  $C_q=\eta_q(H^{1,a}_q(.,L_q))/4$ with large probability. Notice that $C_q>0$ by assumption $[{\cal H}^a]$. We conclude by reminding that we have just proved that $\sqrt{N}(F_N(L_q^N)-F_N(L_{\nu_q^N}))$ is bounded with arbitrarily large probability.
\medskip

Now, for the last term $|L_{\nu_q^N} - L_q|$ of (\ref{truc}), the technique is quite similar. From the first part of the recurrence, taking $\phi=G_{q}$, we have with arbitrarily large probability for $N$ large enough, 
$$\left|\nu_q^N(G_{q})-\alpha\right|=\left|\nu_q^N(G_{q})-\eta_q(G_{q})\right|\leq \frac{C}{\sqrt{N}}.$$
But we also may write
$$\left|\nu_q^N(G_{q})-\alpha\right|=\left|F_N(L_{q})-F_N(L_{\nu_q^N})\right|+o_{L^2}(1/\sqrt{N}).$$
Using the same reasoning as above, we get that for some deterministic constant $C^\prime_q>0$, 
$$
C^\prime_q\left|L_{q}-L_{\nu_q^N}\right|\leq\left|F_N(L_{q})-F_N(L_{\nu_q^N})\right|,
$$ 
and we conclude following the same line. 
\hfill  $\blacksquare$
\medskip

Our last result is then a direct application of Lemma \ref{techlemma} and Proposition \ref{L}.

\begin{cor}\label{equivalence}
For any integer $0\leq q<n$ and for any bounded and measurable function $\phi$, we have
$$\eta_q^{N}(\phi (G_{\eta_q^N}-G_q))=\nu_q^{N}(\phi (G_{\eta_q^N}-G_q))+o_p(1/\sqrt{N}).$$
\end{cor}

\subsection{Proof of Proposition \ref{nec.cond}}\label{preuveH}
We will use the following auxiliary result, which corresponds to Lemma 2.2 in Legoll and Leli\`evre \cite{legoll2010}.

\begin{lem}\label{tony}
Let $f$ denote a mapping from $\R^d$ to $\R$, then the function $F:\R\to\R$ defined by
$$F(L)=\int_{S(x)=L}f(x)\frac{\dbar x}{|DS(x)|}$$
is differentiable with derivative
\begin{eqnarray*}
F'(L)&=&\int_{S(x)=L}\div{\left(f(x)\ \frac{DS(x)}{|DS(x)|^2}\right)}\frac{\dbar x}{|DS(x)|}\\
&=&\int_{S(x)=L}\left[ \frac{DS(x)\cdot Df(x)}{|DS(x)|^2}+f(x)\ \div\left(\frac{DS(x)}{|DS(x)|^2}\right) \right]\frac{\dbar x}{|DS(x)|},
\end{eqnarray*}
provided that the right-hand side is well defined.
\end{lem}

Let us apply this result to the context of Proposition \ref{nec.cond}. We remind the reader that 
\begin{equation}\nonumber
H_q^g(x,L)=\int_{S(x')=L}g(x')K_{q+1}(x,x')\frac{\dbar x'}{|DS(x')|}.
\end{equation}
By the first expression of the derivative in Lemma \ref{tony}, we have 
$$\frac{\partial}{\partial s} H_q^g(x,s)=\int_{S(x')=L} \div_{x'}\left[g(x')\frac{DS(x')K_{q+1}(x,x')}{|DS(x')|^2}\right] \frac{\dbar x'}{|DS(x')|},$$
provided that the right-hand term is well defined. To prove this, notice that 
\begin{align}
\div_{x'}\left[g(x')\frac{DS(x')K_{q+1}(x,x')}{|DS(x')|^2}\right]=&\ g(x')\times\div_{x'}\left[\frac{DS(x')K_{q+1}(x,x')}{|DS(x')|^2}\right]\nonumber\\
&\ + (Dg(x')\cdot DS(x'))\times\frac{K_{q+1}(x,x')}{|DS(x')|^2},\nonumber
\end{align}
where `$\cdot$' stand for the usual scalar product in $\R^d$. For the first term, we use the fact that $g$ is bounded, while for the second one, we apply Cauchy-Schwarz inequality and the inequality between the Euclidean norm $|\cdot|$ and the $L_1$ norm $|\cdot|_1$ to obtain
\begin{align}
&\left|\div_{x'}\left[g(x')\frac{DS(x')K_{q+1}(x,x')}{|DS(x')|^2}\right]\right|\nonumber\\
&\ \leq C\times\left|\div_{x'}\left[\frac{DS(x')K_{q+1}(x,x')}{|DS(x')|^2}\right]\right|+ |\div[g(x')]|_1\times \frac{K_{q+1}(x,x')}{|DS(x')|}.\nonumber
\end{align}  
Concerning the second term, recall that $g$ belongs to 
$${\cal B}_q=\left\{g:\R^d\to\R,\ \exists (g_0\dots,g_{q-1})\in\Ba(\R^d)^q,\ g=K_1(g_0)\cdots K_q(g_{q-1})\right\},$$
so that
\begin{align}
\frac{\partial  g}{\partial x'_j}(x')=\sum_{m=1}^{q}&  K_1(g_0)(x')\dots K_{m-1}(g_{m-2})(x') \left(\int\frac{\partial}{\partial x'_j}K_m(x',x'')g_{m-1}(x'') dx''\right)\nonumber\\
&\  K_{m+1}(g_{m})(x')\dots K_{q}(g_{q-1})(x'),\nonumber
\end{align}
and since all the mappings $g_m$'s are assumed bounded, we get
$$\left|\frac{\partial g}{\partial x'_j}(x')\right|\leq C\sum_{m=1}^{q}\int\left|\frac{\partial}{\partial x'_j}K_m(x',x'')\right|dx'',$$
and finally
$$|\div[g(x')]|_1\leq C\left[\sum_{m=1}^q  \sum_{j=1}^d \int\left|\frac{\partial}{\partial x_j'}K_{m}(x',x'')\right|dx''\right].$$
By the assumption of Proposition \ref{nec.cond}, we deduce that $s\mapsto H_q^g(x,s)$ is differentiable. Moreover, using the mean value theorem, we deduce that
\begin{equation}\nonumber
\left|H_q^g(x,L)-H_q^g(x,L_q)\right|\leq |L-L_q|\times\sup_{s\in(L,L_q)} \left|\frac{\partial}{\partial s} H_q^g(x,s)\right|\leq C|L-L_q|\times h(x),
\end{equation}
with $h\in L^2(\eta)$, so that $[{\cal H}](ii)$ is satisfied.
\hfill  $\blacksquare$ 

\subsection{Proof of the Gaussian case}\label{preuvegauss}
In order to keep the notation as simple as possible, we will explain what happens in dimension $d=2$ only. Thus, the score function is defined, for any $x=(x_1,x_2)\in\R^2$, by $S(x)=x_1/|x|=\cos x$, so that $-1\leq S(x)\leq 1$ and
$$DS(x)=\left[\frac{x_2^2}{|x|^3},\frac{-x_1x_2}{|x|^3}\right]\ \Rightarrow\ |DS(x)|=\frac{|x_2|}{|x|^2}.$$
Hence, denoting $r=|x|$, one has for any $L\in(-1,+1)$
$$S(x)=L\ \Longleftrightarrow\ (x_1,x_2)=(|x|L,\pm |x|\sqrt{1-L^2})=(rL,\pm r\sqrt{1-L^2})$$
so that 
$$|DS(x)|=\frac{\sqrt{1-L^2}}{|x|}\xrightarrow[|x|\to\infty]{}0,$$ 
and, whatever $L$,  $|DS(x)|$ is clearly not bounded from below on the level set $\{S(x)=L\}$. However, for any test function $f$ and any $L\in(-1,+1)$, the coarea formula gives 
\begin{align}
\int_{S(x)=L}f(x)\frac{\dbar x}{|DS(x)|}=&\frac{1}{\sqrt{1-L^2}}\int_0^{\infty}f(rL,r\sqrt{1-L^2})rdr\nonumber\\
&+\frac{1}{\sqrt{1-L^2}}\int_0^{\infty}f(rL,-r\sqrt{1-L^2})rdr.\label{pol}
\end{align}  
In particular, since in this example $X$ is a centered standard Gaussian random vector in $\R^2$, equation (\ref{fy}) shows that the random variable $Y=S(X)$ has density
$$f_Y(s)=\frac{1}{\pi\sqrt{1-s^2}}\ \un_{|s|<1}.$$
This is not surprising since the point $X/|X|$ is uniformly distributed on the unit circle so that $Y=X_1/|X|$ is just the cosine of such a point. Moreover, the transition kernel $K=K_{q+1}$ is a Gaussian transition kernel defined, for the tuning parameter $\s>0$, by
$$K(x,x')=\frac{1+\s^2}{2\pi\s^2}\exp\left(-\frac{1+\s^2}{2\s^2}\left|x'-\frac{x}{\sqrt{1+\s^2}}\right|^2\right).$$ 
Let us recall that point $(i)$ of Assumption [${\cal H}$] requires that
$$\int\eta(dx)\left(\int_{S(x')=L_q}K_{q+1}(x,x')\frac{\dbar x'}{|DS(x')|}\right)^2<\infty.$$
In our context, setting 
$$I=\int_{S(x')=L}K_{q+1}(x,x')\frac{\dbar x'}{|DS(x')|},$$
denoting $\alpha=\sqrt{(1+\s^2)/\s^2}$ and, for any $x=(x_1,x_2)$,
$$A_+=\frac{x_1L+x_2\sqrt{1-L^2}}{\sqrt{1+\s^2}}\hspace{1cm}\mbox{and}\hspace{1cm}A_-=\frac{x_1\sqrt{1-L^2}-x_2L}{\sqrt{1+\s^2}},$$
as well as
$$B_+=\frac{x_1\sqrt{1-L^2}+x_2L}{\sqrt{1+\s^2}}\hspace{1cm}\mbox{and}\hspace{1cm}B_-=\frac{x_1L-x_2\sqrt{1-L^2}}{\sqrt{1+\s^2}},$$
a straightforward computation reveals that
$$I=\frac{\phi(\a A_-)(\phi(\a A_+)+\a A_+\Phi(\a A_+))+\phi(\a B_+)(\phi(\a B_-)+\a B_-\Phi(\a B_-))}{\sqrt{1-L^2}},$$
where $\phi$ and $\Phi$ are respectively the pdf and the cdf of a standard Gaussian random variable. Since $\max(|\a A_+|,|\a B_-|)\leq |x|/\s$, we deduce that   
$$I\leq\frac{2}{\sqrt{1-L^2}}\left(1+\frac{|x|}{\s}\right),$$
and
$$\int\eta(dx)\left(\int_{S(x')=L_q}K_{q+1}(x,x')\frac{\dbar x'}{|DS(x')|}\right)^2\leq \frac{4}{1-L_q^2}\int_{\R^2}\left(1+\frac{|x|}{\s}\right)^2\frac{e^{-\frac{|x|^2}{2}}}{2\pi}dx$$
which is obviously finite, and therefore [${\cal H}$]$(i)$ is satisfied. 
\medskip

In order to prove that [${\cal H}$]$(ii)$ is fulfilled as well, we will make use of Proposition \ref{nec.cond}. Consider first the integral in the sum. From the expression of $K_q$, we have for any $m>0$,
$$
\int\left|\frac{\partial}{\partial x_j'}K_{m}(x',x'')\right|dx''\leq \int C(|x_j'|+|x''_j|)K_m(x',x'')dx''\leq C_1|x'|^{\alpha_1}+C_2,
$$
for $C_1$, $C_2$ and $\alpha_1$ large  enough. Consequently, we have the same type of upper-bound for the whole expression in brackets, meaning that 
$$\sum_{m=1}^q  \sum_{j=1}^d \int\left|\frac{\partial}{\partial x_j'}K_{m}(x',x'')\right|dx''\leq C_1|x'|^{\alpha_1}+C_2.$$
Then, remembering that on the level set $\{S(x')=L\}$, one has $|DS(x')|=\sqrt{1-L^2}/|x'|$, and since $L_q-\delta\leq L\leq L_q+\delta$, we are led to
\begin{eqnarray*}
\int_{S(x')=L} |x'|^{\alpha_1}\frac{K_{q+1}(x,x')}{|DS(x')|^2}\;\dbar x' &\leq& \int_{S(x')=L} K_{q+1}(x,x')\frac{|x'|^{2+\alpha_1}}{1-L^2}\;\dbar x'\\
&\leq&\frac{1}{1-(L_q+\delta)^2}\int_{S(x')=L} K_{q+1}(x_L,x')|x'|^{2+\alpha_1}\;\dbar x'.
\end{eqnarray*}
where $x_L=(|x|L,\pm|x|\sqrt{1-L^2})$ when $x'=(|x'|L,\pm|x'|\sqrt{1-L^2})$. Simple geometric facts indeed show that $|x-x'|\geq|x_L-x'|$ and thus $K_{q+1}(x_L,x')\geq K_{q+1}(x,x')$. Now, by using the same formulation as in (\ref{pol}), the last integral is in fact one dimensional, and is up to a constant a moment of a Gaussian random variable, which is polynomial in its mean :
$$\int_{S(x')=L} |x'|^{\alpha_1}\frac{K_{q+1}(x,x')}{|DS(x')|^2}\;\dbar x' \leq\frac{1}{1-(L_q+\delta)^2} (C_1|x_L|^{\alpha_2}+C_2).$$
Since $|x_L|=|x|$, we have 
$$\int_{S(x')=L} |x'|^{\alpha_1}\frac{K_{q+1}(x,x')}{|DS(x')|^2}\leq\frac{1}{1-(L_q+\delta)^2}  (C_1|x|^{\alpha_2}+C_2),$$
and more generally,
\begin{align}
&\int_{S(x')=L} \left[\sum_{m=1}^q  \sum_{j=1}^d \int\left|\frac{\partial}{\partial x_j'}K_{m}(x',x'')\right|dx''\right]K_{q+1}(x,x')\frac{\dbar x'}{|DS(x')|^2}\nonumber\\
&\ \  \leq \frac{1}{1-(L_q+\delta)^2}  (C_1|x|^{\alpha_2}+C_2).\label{ft}
\end{align}
Hence the second tem in equation (\ref{eq.H2}) is upper bounded by a polynomial in $|x|$, which is of course integrable with respect to the Gaussian measure $\eta$.
\medskip

Now we consider the first term in (\ref{eq.H2}). Observe first that
$$\frac{DS(x')}{|DS(x')|^2}=\left[|x'|,-|x'|\frac{x_1'}{x_2'}\right]\hspace{1cm}\mbox{and}\hspace{1cm}\div\left(\frac{DS(x')}{|DS(x')|^2}\right)=|x'|\frac{x_1'}{(x_2')^2},$$
we get, when $S(x')=L$ and setting $r=|x'|$ as before, 
$$\frac{DS(x')}{|DS(x')|^2}=\left[r,\mp r\frac{L}{\sqrt{1-L^2}}\right]\hspace{1cm}\mbox{and}\hspace{1cm}\div\left(\frac{DS(x')}{|DS(x')|^2}\right)=\frac{L}{1-L^2}.$$
Then
\begin{eqnarray*}
\lefteqn{
\left|\div_{x'}\left[\frac{DS(x')K_{q+1}(x,x')}{|DS(x')|^2}\right] \right|}\\
&\leq& \left|\div\left[\frac{DS(x')}{|DS(x')|^2}\right] \right|K_{q+1}(x,x')
+ \left|\frac{DS(x')}{|DS(x')|^2}\right|\times \left|D_{x'}K_{q+1}(x,x')\right|_1\\
&\leq& \frac{L}{1-L^2}K_{q+1}(x,x')  + \frac{|x'|}{\sqrt{1-L^2}}\times\left|D_{x'}K_{q+1}(x,x')\right|_1.
\end{eqnarray*}
As before, we have
$$\left|D_{x'}K_{q+1}(x,x')\right|_1\leq C(|x|+|x'|)K_{q+1}(x,x'),$$
and, for any $L\in[L_q-\delta,L_q+\delta]$, 
\begin{align}
&\left|\div_{x'}\left[\frac{DS(x')K_{q+1}(x,x')}{|DS(x')|^2}\right] \right|\nonumber\\
&\ \ \leq \left(\frac{L_q+\delta}{1-(L_q+\delta)^2}+ \frac{C|x'|(|x|+|x'|)}{\sqrt{1-(L_q+\delta)^2}}\right)K_{q+1}(x,x').\nonumber
\end{align}
which ensures that, for the term 
$$
\int_{S(x')=L} \left| \div_{x'}\left[\frac{DS(x')K_{q+1}(x,x')}{|DS(x')|^2}\right] \right| \frac{\dbar x'}{|DS(x')|},
$$
we get the same type of upper-bound as in (\ref{ft}). Putting all things together, we have shown inequality (\ref{eq.H2}) of Proposition \ref{nec.cond}, which means that Assumption [${\cal H}$]$(ii)$ is satisfied.
\hfill  $\blacksquare$ 

\subsection{Proof of Lemma \ref{lemmulti}}\label{secmulti}

Our goal is to prove that, for any integer $q$ and any $\ell\in\{1,2,4\}$, we have
$$
N^{\ell/2}~\E\left[\left.\left(U^N_{1-\alpha}-F_N(L_q^N)\right)^\ell\ \right| \ \F_{q-1}^N \right] \xrightarrow[N\to\infty]{\P} 0.
$$

The principle is to sequentially upper-bound the left-hand side. Set $q>0$ and for any $x\in\R^d$ and any $i=1,\dots,N$, let us define the random variables
$$W_{q}^x:=\frac{1}{N}\sum_{j=1}^N \1_{X_q^j=x}\hspace{1cm}\mbox{and for}\ i=1,\dots,N,\hspace{1cm}W_{q}^i:=W_{q}^{X_q^i}.$$
Then, by definition of $F_N(L)$, it is readily seen that
$$\left|U^N_{1-\alpha}-F_N(L_q^N)\right|\leq \sup_{x\in\R^d}W_{q}^x=\max_{1\leq i\leq N}W_{q}^i.$$
First note that, by the assumption on the gradient of $S$, this supremum can only be reached at a sample point $X_{q-1}^i$. Indeed, since the level sets of $S$ have zero Lebesgue measure, then as soon as a transition by the kernel $K$ is accepted, it will give almost surely a unique value of $S$. Hence, the accumulation of the particles $X_q^i$ on a same point $X_{q-1}^j$ can only be caused by resampling.\medskip

Specifically, recall that the multinomial step as described in Section \ref{ad} consists in drawing an $N$-sample $(\hat X_{q-1}^1,\dots,\hat X_{q-1}^N)$ with common distribution
$$\frac{1}{\lceil N\alpha\rceil}\sum_{j:~X_{q-1}^j\geq L^N_{q-1}}~\delta_{X_{q-1}^i}(dx).$$
Let us denote  
$$\{x_1,\dots, x_{\lceil N\alpha\rceil}\}:=\{X_{q-1}^j:\ X_{q-1}^j\geq L^N_{q-1}\}$$
the set of the $\lceil N\alpha\rceil$ particles which are cloned at the multinomial step and, for $1\leq j\leq\lceil N\alpha\rceil$, $N_q^j$ stands for the random number of clones of $x_j$. Said differently, we have 
$$(N_q^1,\dots, N_q^{\lceil N\alpha\rceil})\sim{\cal M}\left(N,\left(\frac{1}{\lceil N\alpha\rceil},\dots,\frac{1}{\lceil N\alpha\rceil}\right)\right),$$
where ${\cal M}(n,(p_1,\dots,p_m))$ is the multinomial law with parameters $n$ and $(p_1,\dots,p_m)$. Then a moment's thought reveals that
$$\max_{1\leq i\leq N}W_{q}^i\leq\max_{1\leq i\leq N}W_{q-1}^i\times\max_{1\leq j\leq\lceil N\alpha\rceil}N_q^j$$
and since the $N_q^j$'s are independent of $\F_{q-1}^N$, we are led to
$$\E\left[\left.\left|U^N_{1-\alpha}-F_N(L_q^N)\right|^\ell\ \right| \ \F_{q-1}^N \right]\leq\left(\max_{1\leq i\leq N}W_{q-1}^i\right)^\ell\E\left[\left(\max_{1\leq j\leq\lceil N\alpha\rceil}N_q^j\right)^\ell\right].$$
Next, Theorem 4.4 in \cite{lucbucket} ensures that
$$\E\left[\left(\max_{1\leq j\leq\lceil N\alpha\rceil}N_q^j\right)^\ell\right]\underset{N\to\infty}{\sim}\left(\frac{\log N}{\log\log N}\right)^\ell.$$
In particular, one has
$$\E\left[\left(\max_{1\leq j\leq\lceil N\alpha\rceil}N_q^j\right)^\ell\right]\leq C_\ell (\log N)^{\ell},$$
and a straightforward induction gives
$$\E\left[\left|U^N_{1-\alpha}-F_N(L_q^N)\right|^\ell\right]\leq C_\ell^q (\log N)^{\ell q}~\E\left[\left(\max_{1\leq i\leq N}W_{0}^i\right)^\ell\right].$$
Finally, as $\eta_0$ is absolutely continuous,  $\max_{1\leq i\leq N}W_{0}^i=1/N$ and we get
$$N^{\ell/2}~\E\left[\left|U^N_{1-\alpha}-F_N(L_q^N)\right|^\ell\right]\leq C_\ell^q (\log N)^{\ell q}N^{-\ell/2}\xrightarrow[N\to\infty]{}0,$$
which concludes the proof of Lemma \ref{lemmulti}.
\hfill  $\blacksquare$

\paragraph{Acknowledgments.} We are greatly indebted to Pierre Del Moral, Fran\c cois Le Gland and Florent Malrieu for valuable comments and insightful suggestions during the redaction of the paper. 

\bibliographystyle{plain}
\bibliography{biblio-cg4}
\end{document}